\documentclass[10pt]{amsart}
\usepackage{amsmath}
\usepackage{amssymb}
\usepackage{amsfonts}
\usepackage{pxfonts}
\usepackage{amsthm}
\usepackage{amsxtra}
\usepackage{fullpage}

\theoremstyle{plain}
\newtheorem{theorem}{Theorem}[section]
\newtheorem{corollary}{Corollary}[section]
\newtheorem{lemma}{Lemma}[section]
\newtheorem{proposition}{Proposition}[section]

\theoremstyle{definition}

\theoremstyle{remark}
\newtheorem{remark}{Remark}[section]

\newcommand{\C}{\mathbb C}

\newcommand{\Z}{\mathbb Z}
\newcommand{\N}{\mathbb N}
\newcommand{\T}{\mathbb T}

\newcommand{\ScSt}{\scriptstyle}
\newcommand{\DySt}{\displaystyle}

\font\teneuf=eufm10  \font\seveneuf=eufm7 \font\fiveeuf=eufm5
\newfam\euffam 
\textfont\euffam=\teneuf \scriptfont\euffam=\seveneuf
\scriptscriptfont\euffam=\fiveeuf

\begin{document}

\title[]{Bi-orthogonal systems on the unit circle, regular semi-classical weights and the discrete Garnier equations}

\author{N.S.~Witte}

\address{Department of Mathematics and Statistics, University of Melbourne, Victoria 3010, Australia}
\email{\tt n.witte@ms.unimelb.edu.au}

\begin{abstract}
We demonstrate that a system of bi-orthogonal polynomials and their associated functions 
corresponding to a regular semi-classical weight on the unit circle constitute a class 
of general classical solutions to the Garnier systems by explicitly constructing its 
Hamiltonian formulation and showing that it coincides with that of a Garnier system. 
Such systems can also be characterised by recurrence relations of the discrete Painlev\'e 
type, for example in the case with one free deformation variable the system was found to 
be characterised by a solution to the discrete fifth Painlev\'e equation. Here we derive 
the canonical forms of the multi-variable generalisation of the discrete fifth Painlev\'e 
equation to the Garnier systems, i.e. for arbitrary numbers of deformation variables.
\end{abstract}

\maketitle

\section{General Structures of Bi-orthogonality}
\setcounter{equation}{0}

Consider a formal complex weight $ w(z) $ and its Fourier coefficients 
$ \{ w_{k} \}_{k \in \Z} $ defined by
\begin{equation}
  w(z) = \sum_{k=-\infty}^{\infty} w_{k}z^k, \quad
  w_{k} = \int_{\mathbb T} \frac{d\zeta}{2\pi i\zeta} w(\zeta)\zeta^{-k} ,
\label{Fcoeff}
\end{equation}
with support on the unit circle where $ \T $ denotes the unit circle $ |\zeta|=1 $ with 
$ \zeta=e^{i\theta} $, $ \theta\in [0,2\pi) $. 
The Toeplitz determinants constructed from the these Fourier coefficients are 
related to averages over the unitary group $ U(n) $ with respect to the Haar measure
by the Heine formula\cite{Ge_1961},\cite{Ge_1977}
\begin{align}
  I_n[w] := \Big \langle \prod_{l=1}^n w(z_l) \Big \rangle_{U(n)}
  & = \frac{1}{(2\pi)^n n!} \int^{2\pi}_{0} d\theta_1 \cdots \int^{2\pi}_{0} d\theta_n 
      \prod_{l=1}^n w(e^{i\theta_l}) \prod_{1\leq j<k\leq n} |e^{i\theta_j}-e^{i\theta_k}|^2
  \nonumber \\
  & = \det[ w_{i-j} ]_{i,j=0,\dots,n-1}, \quad n\geq 1, \quad I_0[w]=1 .
\label{Uavge}
\end{align}
Note that $ \bar{w}_n \neq w_{-n} $ in general and consequently the Toeplitz matrix 
$ ( w_{i-j} )_{i,j=0,\dots,n-1} $ is not necessarily Hermitian.

Notwithstanding the fact that many physically interesting quantities are characterised 
as averages over the unitary group $ U(n) $ we wish to emphasise another perspective.
Intimately connected with such averages are systems of bi-orthogonal polynomials on 
the unit circle which are orthogonal with respect to the weight $ w(z) $ underlying 
the $ U(n) $ average, in the sense of (\ref{orthog:a}) and (\ref{orthog:b}). 
Such systems are equivalent to systems of bi-orthogonal Laurent polynomials, first
introduced by Jones and Thron \cite{JT_1982} and studied subsequently by \cite{JTN_1984}
amongst others, as was shown by Hendriksen and van Rossum \cite{HvR_1986} and Pastro \cite{Pa_1985}
so that all of our conclusions apply equally well to these systems.
A particular class of weights of great interest is the generic or regular semi-classical 
class which are parameterised by the co-ordinates and residues of singular points 
$ \{z_j\}^M_{j=1} $ and $ \{\rho_j\}^M_{j=1} $ respectively (see (\ref{scwgt2}) for the definition).
Some of the relevant properties of this class are summarised in a self-contained way in Section 2. 
The important fact that is relevant here is that 
systems of bi-orthogonal polynomials and their associated functions with such weights
have the property that their monodromy in the complex spectral variable $ z $ is preserved 
under arbitrary deformations of the singularity co-ordinates $ \{z_j\}^M_{j=1} $. 
This fact was first derived in the context of bi-orthogonal polynomial systems on the unit
circle in \cite{FW_2006a}, although it was known to be true for systems of orthogonal 
polynomials on the line due to work by Magnus \cite{Ma_1995}.
This later result was subsequently extended to orthogonal polynomial systems with a certain type of non-generic 
or degenerate semi-classical weight by Bertola, Eynard and Harnad \cite{BEH_2006}. 
The isomonodromic character of bi-orthogonal Laurent polynomial systems with such weights 
has also been established independently by Bertola and Gekhtman \cite{BG_2007}.
Here we specifically demonstrate that bi-orthogonal polynomial systems on the unit circle
with regular semi-classical weights constitute a class of general classical solutions to 
the Garnier systems. 
We carry out this task in Proposition 3.1 by proving that the dynamical equations of the system
with respect to the singularity co-ordinates $ \{z_j\}^M_{j=1} $ are Hamiltonian and
that this coincides precisely with the Garnier system $ {\mathcal G}_N \equiv \{q_j,p_j;K_j,z_j\} $.
The reader should note that the term classical refers to two
distinct notions: semi-classical orthogonal polynomial systems are ones which generalise the
classical systems such as the Hermite, Laguerre, Jacobi or any member of the Askey Table,
whilst a classical solution of a Painlev\'e or Garnier equation is a special solution
constructible from hypergeometric functions and implies a condition on the parameters 
and boundary/initial conditions. 

The significance of this observation is that one can reverse the usual argument and
use the structures derived from approximation theory to deduce new results about
the integrable system. One particular consequence of the identification of the $ U(n) $ 
averages with the Garnier system is their characterisation by recurrence relations of 
the discrete Painlev\'e type. In the case of the simplest $ U(n) $ average with one free 
deformation variable the system was found to be characterised by a solution to the discrete
fifth Painlev\'e equation, see Proposition \ref{discreteG3} in Section 4.
In this section we derive the canonical forms of the higher
analogues of the discrete fifth Painlev\'e for the Garnier systems, i.e. for the
many deformation variable case from the approximation theory structures.
We give the explicit coupled recurrence relations for the two variable Garnier 
system in Proposition \ref{discreteG4}, and the arbitrary variable recurrence relations 
in Proposition \ref{discreteGM}, and these constitute our key results.

We define bi-orthogonal polynomials $ \{\phi_n(z),\bar{\phi}_n(z)\}^{\infty}_{n=0} $
with respect to the weight $ w(z) $ on the unit circle by the orthogonality relation
\begin{equation}
  \int_{\T} \frac{d\zeta}{2\pi i\zeta} w(\zeta)\phi_m(\zeta)\bar{\phi}_n(\bar{\zeta})
   = \delta_{m,n}, \quad m,n \in \Z_{\geq 0} .
\label{onorm}
\end{equation}
Alternatively one can express this definition in terms of orthogonality with 
respect to the monomial basis
\begin{align}
  \int_{\T} \frac{d\zeta}{2\pi i\zeta} w(\zeta)\phi_n(\zeta)\bar{\zeta}^m
  & = \begin{cases}  0 & m<n \\ 1/\kappa_n & m=n \end{cases} ,
\label{orthog:a} \\
  \int_{\T} \frac{d\zeta}{2\pi i\zeta} w(\zeta)\zeta^m\bar{\phi}_n(\bar{\zeta})
  & = \begin{cases}  0 & m<n \\ 1/\kappa_n & m=n \end{cases} .
\label{orthog:b}
\end{align}
Notwithstanding the notation,
$ \bar{\phi}_n $ is not in general equal to the complex conjugate of $ \phi_n $.
The leading coefficients of the polynomials are specified by 
\begin{align}
  \frac{\phi_n(z)}{\kappa_n} &= z^n + \lambda_nz^{n-1} + \mu_nz^{n-2} +\ldots +r_n ,
  \label{phiCff:a} \\
  \frac{\bar{\phi}_n(z)}{\kappa_n} &= z^n + \bar{\lambda}_nz^{n-1} + \bar{\mu}_nz^{n-2} +\ldots +\bar{r}_n , 
  \label{phiCff:b}
\end{align}
where again $ \bar{\lambda}_n $, $ \bar{\mu}_n $, $ \bar{r}_n $ are not in general 
equal to the corresponding complex conjugates of $ \lambda_n $, $ \mu_n $, $ r_n $
respectively. We define a double sequence of $r$-coefficients by
\begin{equation}
  r_n = \frac{\phi_n(0)}{\kappa_n}, \quad 
  \bar{r}_n = \frac{\bar{\phi}_n(0)}{\kappa_n}, \quad n\geq 1, \quad r_0=\bar{r}_0=1 ,
\end{equation}
which differ slightly from the standard definitions of the reflection or Verblunsky 
coefficients $ \alpha_n $, in that $ \alpha_n = -\bar{r}_{n+1} $.
The polynomial coefficients introduced above are related by a system 
of coupled recurrence equations, the first two being \cite{Ge_1961}
\begin{equation}
   \kappa_n^2 - \kappa_{n-1}^2 = \phi_n(0)\bar{\phi}_n(0), \quad
   \lambda_{n}-\lambda_{n-1} = r_{n}\bar{r}_{n-1}.
   \label{l}
\end{equation}

We have a extension of the standard results on the existence of orthogonal polynomial
systems on the unit circle to the bi-orthogonal setting due to Baxter. 
\begin{proposition}[\cite{Ba_1961}]\label{BOPS_exist}
The bi-orthogonal system $ \{\phi_n,\bar{\phi}_n\}^{\infty}_{n=0} $ exists if and only if 
$ I_n \neq 0 $ for all $ n \in \N $.
\end{proposition}

It is a well known result in the theory of Toeplitz determinants \cite{ops_Sz} that
\begin{equation}
   \frac{I_{n+1}[w] I_{n-1}[w]}{(I_{n}[w])^2}
   = 1 - r_{n}\bar{r}_n, \quad n\geq 1 .
\label{I0}
\end{equation}
Rather than dealing with $ \bar{\phi}_n $ we prefer to use the reciprocal polynomial 
$ \phi^*_n(z) $ defined in terms of the $ n$th degree polynomial $ \bar{\phi}_n(z) $ by
\begin{equation}
  \phi^*_n(z) := z^n\bar{\phi}_n(1/z) .
\label{recip}
\end{equation}
The generating function of the Toeplitz elements, known as the Carath\'eodory function
\begin{equation}
   F(z) \coloneqq \int_{\T}\frac{d\zeta}{2\pi i\zeta}\frac{\zeta+z}{\zeta-z}w(\zeta) ,
\label{Cfun}
\end{equation}
will also feature prominently in our work.
The fundamental object identified in the study \cite{FW_2006a} is the $ 2\times 2$ matrix
\begin{equation}
   Y_n(z;t) := 
   \begin{pmatrix}
          \phi_n(z)   &  \epsilon_n(z)/w(z) \cr
          \phi^*_n(z) & -\epsilon^*_n(z)/w(z) \cr
   \end{pmatrix} ,
\label{Ydefn}
\end{equation}
where the associated functions or functions of the second kind are defined by
\begin{align}
   \epsilon_n(z)
   &:= \int_{\T}\frac{d\zeta}{2\pi i\zeta}\frac{\zeta+z}{\zeta-z}w(\zeta)
                   \phi_n(\zeta), \quad n\geq 1, \quad \epsilon_0(z) = \kappa_0[w_0+F(z)] ,
   \label{eps:a} \\
   \epsilon^*_n(z)
   &:= -z^n\int_{\T}\frac{d\zeta}{2\pi i\zeta}\frac{\zeta+z}{\zeta-z}w(\zeta)
                   \bar{\phi}_n(\bar{\zeta}), \quad n\geq 1, \quad \epsilon^*_0(z) = \kappa_0[w_0-F(z)] .
   \label{eps:b}
\end{align}

Solutions to the orthogonality relations yield the following determinantal and integral
representations for the polynomials,
\begin{align}
  \phi_{n}(z) & = \frac{\kappa_n}{I_{n}}
         \det \begin{pmatrix}
                     w_{0}      & \ldots & w_{-j}       & \ldots        & w_{-n} \cr
                     \vdots     & \vdots & \vdots & \vdots      & \vdots \cr
                     w_{n-1}    & \ldots & w_{n-j-1}    & \ldots        & w_{-1} \cr
                     1  & \ldots & z^j  & \ldots        & z^n  \cr
              \end{pmatrix}
    = (-1)^n\kappa_n\frac{I_{n}[w(\zeta)(\zeta-z)]}{I_{n}[w(\zeta)]} ,
\label{phiRep} \\
  \phi^*_{n}(z) & = \frac{\kappa_n}{I_{n}}
         \det \begin{pmatrix}
                     w_{0}      & \ldots & w_{-n+1}     & z^n  \cr
                     \vdots     & \vdots & \vdots       & \vdots \cr
                     w_{n-j}    & \ldots & w_{-j+1}     & z^j  \cr
                     \vdots     & \vdots & \vdots       & \vdots \cr
                     w_{n}      & \ldots & w_{1}        & 1  \cr
              \end{pmatrix}
    = \kappa_n\frac{I_{n}[w(\zeta)(1-z\zeta^{-1})]}{I_{n}[w(\zeta)]} .
\label{phiSRep}
\end{align}
The associated functions have representations analogous to (\ref{phiRep},\ref{phiSRep})
\begin{align}
  \frac{\kappa_n}{2}\epsilon_{n}(z)
 & = \frac{1}{2I_{n+1}}
          \det \begin{pmatrix}
                     w_{0}      & \ldots & w_{-j}       & \ldots        & w_{-n} \cr
                     \vdots     & \vdots & \vdots & \vdots      & \vdots \cr
                     w_{n-1}    & \ldots & w_{n-j-1}    & \ldots        & w_{-1} \cr
                     g_0  & \ldots & g_j  & \ldots        & g_n  \cr
              \end{pmatrix}
    = z^n\frac{I_{n+1}[w(\zeta)(1-z\zeta^{-1})^{-1}]}{I_{n+1}[w(\zeta)]} ,
\label{epsRep} \\
  \frac{\kappa_n}{2}\epsilon^*_{n}(z)
 & = (-1)^{n+1}\frac{1}{2I_{n+1}}
         \det \begin{pmatrix}
                     w_{0}      & \ldots & w_{-n+1}     & g_n  \cr
                     \vdots     & \vdots & \vdots       & \vdots \cr
                     w_{n-j}    & \ldots & w_{-j+1}     & g_j  \cr
                     \vdots     & \vdots & \vdots       & \vdots \cr
                     w_{n}      & \ldots & w_{1}        & g_0  \cr
              \end{pmatrix}
    = (-z)^{n+1}\frac{I_{n+1}[w(\zeta)(\zeta-z)^{-1}]}{I_{n+1}[w(\zeta)]} .
\label{epsSRep}
\end{align}
where
\begin{equation}
	g_n(z) := 2z \int_{\T}\frac{d\zeta}{2\pi i\zeta}w(\zeta)\frac{\zeta^n}{\zeta-z}, \quad n\geq 0 ,
\end{equation}
for $ |z| \neq 1 $.

From the general properties of bi-orthogonality we can deduce that the matrix $ Y_n $ 
obeys a difference system.
\begin{proposition}[\cite{Ge_1961}]
Assuming $ \kappa_n \neq 0 $ (or equivalently $ I_n \neq 0 $) for $ n\in \N $
the matrix $ Y_n $ satisfies the recurrence relation in $ n $
\begin{equation}
   Y_{n+1} := K_n Y_{n}
   = \frac{1}{\kappa_n}
       \begin{pmatrix}
              \kappa_{n+1} z   & \phi_{n+1}(0) \cr
              \bar{\phi}_{n+1}(0) z & \kappa_{n+1} \cr
       \end{pmatrix} Y_{n} .
\label{Yrecur}
\end{equation}
\end{proposition}

\begin{theorem}[\cite{Ge_1961}]
The Casoratians of the solutions $ \phi_n, \phi^*_n, \epsilon_n, \epsilon^*_{n} $ 
to the above recurrence relations are
\begin{align}
  \phi_{n+1}(z)\epsilon_n(z) - \epsilon_{n+1}(z)\phi_n(z)
  & = 2\frac{\phi_{n+1}(0)}{\kappa_n}z^n ,
  \label{Cas:a} \\
  \phi^*_{n+1}(z)\epsilon^*_n(z) - \epsilon^*_{n+1}(z)\phi^*_n(z)
  & = 2\frac{\bar{\phi}_{n+1}(0)}{\kappa_n}z^{n+1} ,
  \label{Cas:b} \\
  \phi_{n}(z)\epsilon^*_n(z) + \epsilon_{n}(z)\phi^*_n(z)
  & = 2z^n .
  \label{Cas:c}
\end{align}
\end{theorem}

Under quite general conditions the matrix system $ Y_n $ obeys the following spectral differential system. 
\begin{proposition}[\cite{IW_2001},\cite{FW_2006a}]\label{BOPS_specDer}
Assume that the weight satisfies the moment conditions
\begin{equation}
	\int_{\T} \frac{d\zeta}{2\pi i\zeta} w(\zeta)
  \frac{\frac{\DySt d}{\DySt dz}\log w(z)-\frac{\DySt d}{\DySt d\zeta}\log w(\zeta)}{z-\zeta}\zeta^j
  \neq \infty, \quad j \in \Z .
\end{equation}
Then the matrix $ Y_n $ satisfies the differential relation with respect to the
spectral variable $ z $ 
\begin{equation}
   \frac{d}{dz}Y_{n} := A_n Y_{n}
   = \frac{1}{W(z)}
       \begin{pmatrix}
              -\left[ \Omega_n(z)+V(z)
                     -\dfrac{\kappa_{n+1}}{\kappa_n}z\Theta_n(z)
               \right]
            & \dfrac{\phi_{n+1}(0)}{\kappa_n}\Theta_n(z)
            \cr
              -\dfrac{\bar{\phi}_{n+1}(0)}{\kappa_n}z\Theta^*_n(z)
            &  \Omega^*_n(z)-V(z)
                     -\dfrac{\kappa_{n+1}}{\kappa_n}\Theta^*_n(z)
            \cr
       \end{pmatrix} Y_{n} .
\label{YspecDer}
\end{equation}
\end{proposition}
The utility of such a parameterisation of the spectral matrix $ A_n $ will be 
evident when we make the specialisation to the regular semi-classical weights.
We will refer to $ \Theta_n, \Omega_n, \Theta^*_n, \Omega^*_n $ as spectral
coefficients.

The scalar differential equation system corresponding to the above matrix system
is specified in the following result.
\begin{proposition}\label{scalarODE}
The components of the matrix $ Y_n $ satisfy two second-order scalar 
ordinary differential equations in the spectral variable: $ \phi_n(z) $ or
$ \epsilon_n(z)/w(z) $ satisfy 
\begin{equation}
	\phi_n''+p_1\phi_n'+p_2\phi_n = 0 ,
\label{2ODE:a}
\end{equation}
while $ \phi^*_n(z) $ and $ -\epsilon^*_n(z)/w(z) $ satisfy
\begin{equation}
	\phi^{*''}_n+p^*_1\phi^{*'}_n+p^*_2\phi^*_n = 0 .
\label{2ODE:b}
\end{equation}
The coefficients of the scalar second-order differential equations are 
\begin{equation}
	p_1(z) = \frac{W'}{W}-\frac{\Theta_n'}{\Theta_n}+\frac{2V}{W}-\frac{n}{z} ,
\label{2ODEcff:a}
\end{equation}
and
\begin{multline}
	p_2(z) = \frac{\Theta_n(\Omega_n'+V')-\Theta_n'(\Omega_n+V)}{W\Theta_n}
	 -\frac{\kappa_{n+1}}{\kappa_n}\frac{\Theta_n}{W}
	\\
	 -\frac{\left[\Omega_n+V-\frac{\DySt \kappa_{n+1}}{\DySt \kappa_n}z\Theta_n\right]
	        \left[\Omega^*_n-V-\frac{\DySt \kappa_{n+1}}{\DySt \kappa_n}\Theta^*_n\right]
	       }{W^2}
	 +\frac{\phi_{n+1}(0)\bar{\phi}_{n+1}(0)z\Theta_n\Theta^*_n}{\kappa^2_nW^2} ,
\label{2ODEcff:b}
\end{multline}
and
\begin{equation}
	p^*_1(z) = \frac{W'}{W}-\frac{\Theta^{*'}_n}{\Theta^*_n}+\frac{2V}{W}-\frac{n+1}{z} ,
\label{2ODEcff:c}
\end{equation}
and 
\begin{multline}
	p^*_2(z) = \frac{(z^{-1}\Theta^*_n+\Theta^{*'}_n)(\Omega^*_n-V)-\Theta^*_n(\Omega^{*'}_n-V')}{W\Theta^*_n}
	 -\frac{\kappa_{n+1}}{\kappa_n}\frac{\Theta^*_n}{zW}
	\\
	 -\frac{\left[\Omega_n+V-\frac{\DySt \kappa_{n+1}}{\DySt \kappa_n}z\Theta_n\right]
	        \left[\Omega^*_n-V-\frac{\DySt \kappa_{n+1}}{\DySt \kappa_n}\Theta^*_n\right]
	       }{W^2}
	 +\frac{\phi_{n+1}(0)\bar{\phi}_{n+1}(0)z\Theta_n\Theta^*_n}{\kappa^2_nW^2} .
\label{2ODEcff:d}
\end{multline}
\end{proposition}
\begin{proof}
The two ordinary differential equations (\ref{2ODE:a}) and (\ref{2ODE:b}) follow 
from the elimination of $ \phi^*_n $ and $ \phi_n $ in (\ref{YspecDer}) respectively. 
\end{proof}

A consequence of the compatibility between the differential relations (\ref{YspecDer}) and 
the recurrence relations (\ref{Yrecur}) is the following collection of recurrence relations for 
the spectral coefficients $ \Omega_n, \Omega^*_n, \Theta_n, \Theta^*_n $.
\begin{proposition}[\cite{FW_2006a}]\label{Linear1}
Given the conditions of Proposition \ref{BOPS_exist} the spectral coefficients
$ \{ \Omega_n(z), \Omega^*_n(z), \Theta_n(z), \Theta^*_n(z) \}^{\infty}_{n=0} $
satisfy the following recurrence relations in $ n $
\begin{equation}
  \Omega_n(z) + \Omega_{n-1}(z)
  - \left( \frac{\phi_{n+1}(0)}{\phi_{n}(0)}+\frac{\kappa_{n+1}}{\kappa_{n}}z
    \right)\Theta_n(z) + (n-1)\frac{W(z)}{z} = 0 ,
\label{rrCf:a}
\end{equation}
\begin{equation}
   \left( \frac{\phi_{n+1}(0)}{\phi_{n}(0)}+\frac{\kappa_{n+1}}{\kappa_{n}}z
   \right) (\Omega_{n-1}(z) - \Omega_{n}(z))
   + \frac{\kappa_{n}\phi_{n+2}(0)}{\kappa_{n+1}\phi_{n+1}(0)}z\Theta_{n+1}(z)
   - \frac{\kappa_{n-1}\phi_{n+1}(0)}{\kappa_{n}\phi_{n}(0)}z\Theta_{n-1}(z)
   - \frac{\phi_{n+1}(0)}{\phi_{n}(0)}\frac{W(z)}{z} = 0 ,
\label{rrCf:b}
\end{equation}
\begin{equation}
  \Omega^*_n(z) + \Omega^*_{n-1}(z)
  - \left( \frac{\kappa_{n+1}}{\kappa_{n}}+\frac{\bar{\phi}_{n+1}(0)}{\bar{\phi}_{n}(0)}z
    \right)\Theta^*_n(z) - n\frac{W(z)}{z} = 0 ,
\label{rrCf:c}
\end{equation}
\begin{equation}
   \left( \frac{\kappa_{n+1}}{\kappa_{n}}+\frac{\bar{\phi}_{n+1}(0)}{\bar{\phi}_{n}(0)}z
   \right) (\Omega^*_{n-1}(z) - \Omega^*_{n}(z))
   + \frac{\kappa_{n}\bar{\phi}_{n+2}(0)}{\kappa_{n+1}\bar{\phi}_{n+1}(0)}
      z\Theta^*_{n+1}(z)
   - \frac{\kappa_{n-1}\bar{\phi}_{n+1}(0)}{\kappa_{n}\bar{\phi}_{n}(0)}
      z\Theta^*_{n-1}(z)
   + \frac{\kappa_{n+1}}{\kappa_{n}}\frac{W(z)}{z} = 0 ,
\label{rrCf:d}
\end{equation}
\begin{equation}
  \Omega_{n+1}(z) + \Omega^*_{n}(z)
  - \left( \frac{\phi_{n+2}(0)}{\phi_{n+1}(0)}+\frac{\kappa_{n+2}}{\kappa_{n+1}}z
    \right)\Theta_{n+1}(z)
  + \frac{\kappa_{n+1}}{\kappa_{n}}(z\Theta_n(z)-\Theta^*_n(z)) = 0 ,
\label{rrCf:e}
\end{equation}
\begin{equation}
   \Omega_{n}(z) - \Omega_{n+1}(z)
   + \frac{\kappa_{n+2}}{\kappa_{n+1}}
    \left( z+\frac{\bar{\phi}_{n+1}(0)}{\kappa_{n+1}}\frac{\phi_{n+2}(0)}{\kappa_{n+2}}
    \right)\Theta_{n+1}(z)
   + \frac{\phi_{n+1}(0)\bar{\phi}_{n+1}(0)}{\kappa_{n+1}\kappa_{n}}\Theta^*_n(z)
   - \frac{\kappa_{n+1}}{\kappa_{n}}z\Theta_{n}(z)
   - \frac{W(z)}{z} = 0 ,
\label{rrCf:f}
\end{equation}
\begin{equation}
  \Omega^*_{n+1}(z) + \Omega_{n}(z)
  - \left( \frac{\kappa_{n+2}}{\kappa_{n+1}}
           +\frac{\bar{\phi}_{n+2}(0)}{\bar{\phi}_{n+1}(0)}z
    \right)\Theta^*_{n+1}(z)
  - \frac{\kappa_{n+1}}{\kappa_{n}}(z\Theta_n(z)-\Theta^*_n(z)) - \frac{W(z)}{z}= 0 ,
\label{rrCf:g}
\end{equation}
\begin{equation}
   \Omega^*_{n}(z) - \Omega^*_{n+1}(z)
   + \frac{\kappa_{n+2}}{\kappa_{n+1}}
    \left( 1+\frac{\phi_{n+1}(0)}{\kappa_{n+1}}\frac{\bar{\phi}_{n+2}(0)}{\kappa_{n+2}}z
    \right)\Theta^*_{n+1}(z)
   + \frac{\phi_{n+1}(0)\bar{\phi}_{n+1}(0)}{\kappa_{n+1}\kappa_{n}}z\Theta_n(z)
   - \frac{\kappa_{n+1}}{\kappa_{n}}\Theta^*_{n}(z) = 0 .
\label{rrCf:h}
\end{equation}
\end{proposition}

The spectral coefficients $ \Theta_{n},\Omega_{n} $ are related to their "conjugates"
$ \Theta^*_{n},\Omega^*_{n} $ via a number of functional (or recurrence) relations so that
one can characterise the system in terms of either set. We will refer to these as transition
relations.
\begin{corollary}[\cite{FW_2006a}]\label{Linear2}
The spectral coefficients are inter-related through the following equations
\begin{gather}
   \frac{\bar{\phi}_{n+1}(0)}{\bar{\phi}_{n}(0)}z\Theta^*_n(z)
     - \frac{\kappa_{n}}{\kappa_{n-1}} \Theta^*_{n-1}(z)
   = \frac{\phi_{n+1}(0)}{\phi_{n}(0)}\Theta_{n}(z)
     - \frac{\kappa_{n}}{\kappa_{n-1}}z\Theta_{n-1}(z) ,
\label{rrCf:i} \\
   \Omega^*_{n}(z)-\frac{\kappa_{n+1}}{\kappa_{n}}\Theta^*_n(z)
    = \Omega_{n}(z)-\frac{\kappa_{n+1}}{\kappa_{n}}z\Theta_{n}(z)+n\frac{W(z)}{z} ,
\label{rrCf:j} \\
   \Omega^*_{n}(z)+\Omega_{n}(z)
    = \frac{\kappa^2_{n}}{\kappa^2_{n+1}}
      \left[\frac{\phi_{n+2}(0)}{\phi_{n+1}(0)}\Theta_{n+1}(z)
             + \frac{\kappa_{n+1}}{\kappa_{n}} \Theta^*_{n}(z)\right]+\frac{W(z)}{z} .
\label{rrCf:k}
\end{gather}
\end{corollary}

We will require the leading order terms in expansions of
$ \phi_n(z), \phi^*_n(z), \epsilon_n(z), \epsilon^*_{n}(z) $ both inside and
outside the unit circle. The following Corollary extends the preliminary results
reported in \cite{FW_2006a}.
\begin{corollary}[\cite{FW_2006a}]
The bi-orthogonal polynomials $ \phi_n(z), \phi^*_n(z) $ have the following
expansions for $ |z| < \delta^{-} < 1 $
\begin{equation}
   \dfrac{1}{\kappa_n}\phi_n(z) =
       r_n + \left[r_{n-1}+r_n\bar{\lambda}_{n-1}\right]z
       + \left[r_{n-2}+r_{n-1}\bar{\lambda}_{n-2}
               +r_{n}\bar{\mu}_{n-1} \right]z^2
       + {\rm O}(z^3) ,
\label{phiexp:a}
\end{equation}
\begin{equation}
   \dfrac{1}{\kappa_n}\phi^*_n(z) =
      1 + \bar{\lambda}_n z + \bar{\mu}_n z^2 + \bar{\nu}_n z^3 + {\rm O}(z^{3}) ,
\label{phiexp:b}
\end{equation}
whilst the associated functions have the expansions
\begin{multline}
   \dfrac{\kappa_n}{2}\epsilon_n(z) =
      z^n - \bar{\lambda}_{n+1}z^{n+1}
   + \left[\bar{\lambda}_{n+1}\bar{\lambda}_{n+2}-\bar{\mu}_{n+2}\right]z^{n+2} \\
   + \left[\bar{\lambda}_{n+1}\bar{\mu}_{n+3}+\bar{\lambda}_{n+3}\bar{\mu}_{n+2}
           -\bar{\nu}_{n+3}-\bar{\lambda}_{n+1}\bar{\lambda}_{n+2}\bar{\lambda}_{n+3} \right]z^{n+3} 
       + {\rm O}(z^{n+4}) ,
\label{epsexp:a}
\end{multline}
\begin{equation}
   \dfrac{\kappa_n}{2}\epsilon^*_n(z) = \bar{r}_{n+1}z^{n+1}
       + \left[ \bar{r}_{n+2}-\bar{r}_{n+1}\bar{\lambda}_{n+2} \right]z^{n+2}
       + \left[ \bar{r}_{n+3}-\bar{r}_{n+2}\bar{\lambda}_{n+3}-\bar{r}_{n+1}\bar{\mu}_{n+3}
               +\bar{r}_{n+1}\bar{\lambda}_{n+2}\bar{\lambda}_{n+3}
         \right]z^{n+3}
       + {\rm O}(z^{n+4}) .
\label{epsexp:b}
\end{equation}
The large argument expansions $ |z| > \delta^{+} > 1 $ of $ \phi_n(z), \phi^*_n(z) $ are 
\begin{equation}
   \dfrac{1}{\kappa_n}\phi_n(z) =
      z^n + \lambda_n z^{n-1} + \mu_n z^{n-2} + \nu_n z^{n-3} + {\rm O}(z^{n-3}) ,
\label{phiexp:c}
\end{equation}
\begin{equation}
   \dfrac{1}{\kappa_n}\phi^*_n(z) =
      \bar{r}_n z^n
       + \left[\bar{r}_{n-1}+\bar{r}_n\lambda_{n-1}\right]z^{n-1}
       + \left[\bar{r}_{n-2}+\bar{r}_{n-1}\lambda_{n-2}+\bar{r}_{n}\mu_{n-1}\right]z^{n-2} 
       + {\rm O}(z^{n-3}) ,
\label{phiexp:d}
\end{equation}
whilst the associated functions have the expansions
\begin{equation}
   \dfrac{\kappa_n}{2}\epsilon_n(z) = r_{n+1}z^{-1}
       + \left[r_{n+2}-r_{n+1}\lambda_{n+2}\right]z^{-2}
       + \left[r_{n+3}-r_{n+2}\lambda_{n+3}-r_{n+1}\mu_{n+3}
                +r_{n+1}\lambda_{n+2}\lambda_{n+3}\right]z^{-3}
       + {\rm O}(z^{-4}) ,
\label{epsexp:c}
\end{equation}
\begin{equation}
   \dfrac{\kappa_n}{2}\epsilon^*_n(z) =
       1 - \lambda_{n+1}z^{-1} + \left[\lambda_{n+2}\lambda_{n+1}-\mu_{n+2}\right]z^{-2}
   + \left[\lambda_{n+1}\mu_{n+3}+\lambda_{n+3}\mu_{n+2}
           -\nu_{n+3}-\lambda_{n+1}\lambda_{n+2}\lambda_{n+3}\right]z^{-3}       
        + {\rm O}(z^{-4}) .
\label{epsexp:d}
\end{equation}
\end{corollary}
\begin{proof}
Expansions (\ref{phiexp:a}) and (\ref{phiexp:d}) are found by differentiating 
\begin{equation}
	\kappa_n\phi_{n+1}(z) = \kappa_{n+1}z\phi_n(z)+\phi_{n+1}(0)\phi^*_n(z) ,
\end{equation}
and
\begin{equation}
	\phi_{n+1}(0)\bar{\phi}_{n+1}(z) = \kappa_{n+1}z^{n+1}\phi_{n+1}(z^{-1})-\kappa_{n}z^n\phi_n(z^{-1}) ,
\end{equation}
repeatedly, respectively, and setting the argument to zero.
Expansion (\ref{epsexp:a}) can be found using 
\begin{equation}
  z^n = \frac{\bar{\phi}_n(z)}{\kappa_n} - \bar{\lambda}_{n}\frac{\bar{\phi}_{n-1}(z)}{\kappa_{n-1}}
   + \left[\bar{\lambda}_{n}\bar{\lambda}_{n-1}-\bar{\mu}_{n}\right]\frac{\bar{\phi}_{n-2}(z)}{\kappa_{n-2}}
   + \left[\bar{\mu}_{n}\bar{\lambda}_{n-2}+\bar{\mu}_{n-1}\bar{\lambda}_{n}
           -\bar{\nu}_{n}-\bar{\lambda}_{n}\bar{\lambda}_{n-1}\bar{\lambda}_{n-2}\right]
             \frac{\bar{\phi}_{n-3}(z)}{\kappa_{n-3}} 
   + \Pi_{n-4} .
\end {equation}
Expansion (\ref{epsexp:c}) is found by making use of
\begin{equation}
	z\phi_n(z) = \frac{\kappa_n}{\kappa_{n+1}}\phi_{n+1}(z)
	   -\frac{\phi_{n+1}(0)}{\kappa_n\kappa_{n+1}}\sum^{n}_{j=0}\bar{\phi}_j(0)\phi_j(z) ,
\end{equation}
repeatedly. Expansion (\ref{epsexp:b}) can be found using the conjugate 
analogue of the above equation.
\end{proof}

\section{The Regular Semi-classical Class of Weights}
\setcounter{equation}{0}
Of direct relevance to integrable systems is the regular semi-classical class, 
characterised by a special structure of their logarithmic derivatives
\begin{equation}
  \frac{1}{w(z)}\frac{d}{dz}w(z) = \frac{2V(z)}{W(z)}
  = \sum^M_{j=1}\frac{\rho_j}{z-z_j}, \quad \rho_j \in \C,
\label{scwgt2}
\end{equation}
and its degenerate cases.
Here $ V(z) $, $ W(z) $ are polynomials satisfying the following generic conditions 
for the regular semi-classical class -
\begin{enumerate}
 \item[(i)]
  $ {\rm deg}\;(W) \geq 2 $,
 \item[(ii)]
  $ {\rm deg}\;(V) < {\rm deg}\;(W)=M $,
 \item[(iii)]
  the $ M $ zeros of $ W(z) $, $ \{z_1, \ldots ,z_M\} $ are distinct, and
 \item[(iv)]
  the residues $ \rho_j = 2V(z_j)/W'(z_j) \notin \Z_{\geq 0} $.
\end{enumerate} 
We have the expansion of the denominator in terms of elementary symmetric functions
\begin{equation}
  W(z) = \prod^M_{j=1}(z-z_j) = \sum^{M}_{l=0} (-)^le_l[z_1,\dots,z_M] z^{M-l}, \quad e_0=1 ,
\end{equation}
and of the numerator
\begin{equation}
  2V(z) = \sum^{M-1}_{l=0} (-)^lm_l[z_1,\dots,z_M] z^{M-1-l}, \quad m_0 = \sum^M_{j=1} \rho_j .
\end{equation}
One explicit example, however not the most general form, of such a weight is
the generalised Jacobi weight
\begin{equation}
  w(z) = \prod^M_{j=1}(z-z_j)^{\rho_j}, \quad \rho_j \in \C, \quad supp(wdz)=\T .
\label{gJwgt}
\end{equation}

\begin{lemma}[\cite{La_1972b},\cite{Ma_1995a},\cite{FW_2006a}]
Let the weight $ w(z) $ satisfy the conditions of Proposition \ref{BOPS_specDer} and
$ w(e^{2\pi i})=w(1) $. The Carath\'eodory function (\ref{Cfun}) satisfies the first order linear
ordinary differential equation  
\begin{equation}
   W(z)F'(z) = 2V(z)F(z)+U(z) ,
\label{FD}                  
\end{equation}
where $ U(z) $ is a polynomial in $ z $. 
\end{lemma}

Note that we do not assume one of the singularities is located at the origin and the next 
result is a variant of Proposition 3.1 in \cite{FW_2006a}, which did make that assumption.
\begin{proposition}[\cite{FW_2006a}]
For regular semi-classical weights (\ref{gJwgt}), the functions $ z\Theta_n(z) $,
$ z\Theta^*_n(z) $, $ z\Omega_n(z) $ and $ z\Omega^*_n(z) $ in (\ref{YspecDer}) are
polynomials of degree $ {\rm deg}\; z\Omega_n(z)={\rm deg}\; z\Omega^*_n(z)=M $,
$ {\rm deg}\; z\Theta_n(z)={\rm deg}\; z\Theta^*_n(z)=M-1 $, independent of $ n $.
\end{proposition}

Because of the assumption $ z_j\neq 0 $ the following result also differs in detail with 
the corresponding result in \cite{FW_2006a}.
\begin{proposition}[\cite{FW_2006a}]
The spectral coefficients have terminating expansions in the interior domain of the unit circle
about $ z=0 $ with the explicit forms
\begin{equation}
 (-1)^M\frac{\phi_{n+1}(0)}{\phi_{n}(0)}\Theta_n(z) =
    -ne_{M}z^{-1}
    + \bigg\{ ne_{M-1}-m_{M-1}
       + e_{M}\left[(n+1)\bar{\lambda}_{n+1}-(n-1)\left(\bar{\lambda}_{n-1}+\dfrac{r_{n-1}}{r_{n}}\right) \right]
      \bigg\} + {\rm O}(z) ,
\label{Thexp:b}
\end{equation}
\begin{equation}
  (-1)^M\Omega_n(z) = -ne_{M}z^{-1}
  + \bigg\{ ne_{M-1}-\frac{1}{2}m_{M-1}
           +e_{M}\left[(n+1)\bar{\lambda}_{n+1}-n\left(\bar{\lambda}_{n}+\frac{r_n}{r_{n+1}}\right) \right]
   \bigg\} + {\rm O}(z) ,
\label{Omexp:b}
\end{equation}
\begin{equation}
  (-1)^M\frac{\kappa_{n+1}}{\kappa_n}\Theta^*_n(z) = (n+1)e_{M}z^{-1}
 + \bigg\{-(n+1)e_{M-1}+m_{M-1}+e_{M}\left[ (n+2)\left(\frac{\bar{r}_{n+2}}{\bar{r}_{n+1}}-\bar{\lambda}_{n+2}\right)+n\bar{\lambda}_n \right]
   \bigg\} + {\rm O}(z) ,
\label{ThSexp:b}
\end{equation}
\begin{equation}
  (-1)^M\Omega^*_n(z) = (n+1)e_{M}z^{-1}
 + \bigg\{ \frac{1}{2}m_{M-1}-(n+1)e_{M-1}+e_{M}\left[ (n+1)\bar{\lambda}_{n+1}+(n+2)\left(\frac{\bar{r}_{n+2}}{\bar{r}_{n+1}}-\bar{\lambda}_{n+2}\right) \right]
   \bigg\}
 + {\rm O}(z) ,
\label{OmSexp:b}
\end{equation}
and in the exterior domain of the unit circle about $ z=\infty $ with the explicit forms
\begin{equation}
  \frac{\kappa_{n+1}}{\kappa_n}\Theta_n(z) =
    (n+1+m_0)z^{M-2}
    + \bigg\{ -(n+1)e_1-m_1+(n+2+m_0)\left[\frac{r_{n+2}}{r_{n+1}}-\lambda_{n+2}\right]
           +(n+m_0)\lambda_{n} \bigg\}z^{M-3} + {\rm O}(z^{M-4}) ,
\label{Thexp:a}
\end{equation}
\begin{equation}
  \Omega_n(z) =
 (1+\frac{1}{2} m_0)z^{M-1}
 + \bigg\{ -e_{1}-\frac{1}{2}m_{1}+(n+1+m_0)\lambda_{n+1}-(n+2+m_0)\left[\lambda_{n+2}-\frac{r_{n+2}}{r_{n+1}} \right]
   \bigg\}z^{M-2} + {\rm O}(z^{M-3}) ,
\label{Omexp:a}
\end{equation}
\begin{equation}
  \frac{\bar{\phi}_{n+1}(0)}{\bar{\phi}_{n}(0)}\Theta^*_n(z) = -(n+m_0)z^{M-2}
 + \bigg\{ ne_{1}+m_{1}+(n+1+m_0)\lambda_{n+1}-(n-1+m_0)\left[\lambda_{n-1}+\frac{\bar{r}_{n-1}}{\bar{r}_n}\right] 
   \bigg\}z^{M-3} + {\rm O}(z^{M-4}) ,
\label{ThSexp:a}
\end{equation}
\begin{equation}
  \Omega^*_n(z) = -\frac{1}{2}m_0 z^{M-1}
 + \bigg\{ \frac{1}{2}m_1+(n+1+m_0)\lambda_{n+1}-(n+m_0)\left[ \lambda_{n}+\frac{\bar{r}_n}{\bar{r}_{n+1}} \right]
   \bigg\}z^{M-2}
 + {\rm O}(z^{M-3}) .
\label{OmSexp:a}
\end{equation}
\end{proposition}
\begin{proof}
These expansions follow from the inversion of (\ref{YspecDer}), namely the formulae
\begin{align}
  2\frac{\phi_{n+1}(0)}{\kappa_n}z^n\Theta_n
  & = W\left[-\epsilon'_n\phi_n+\epsilon_n\phi'_n \right]+2V\epsilon_n\phi_n ,
\\
  2\frac{\phi_{n+1}(0)}{\kappa_n}z^n\Omega_n
  & = W\left[-\epsilon'_n\phi_{n+1}+\epsilon_{n+1}\phi'_n \right]
      +V\left[ \epsilon_{n+1}\phi_n+\epsilon_n\phi_{n+1} \right] ,
\\
  2\frac{\bar{\phi}_{n+1}(0)}{\kappa_n}z^{n+1}\Theta^*_n
  & = W\left[ \epsilon^{*'}_n\phi^*_n-\epsilon^*_n\phi^{*'}_n \right]-2V\epsilon^*_n\phi^*_n ,
\\
  2\frac{\bar{\phi}_{n+1}(0)}{\kappa_n}z^{n+1}\Omega^*_n
  & = W\left[ \epsilon^{*'}_n\phi^*_{n+1}-\epsilon^*_{n+1}\phi^{*'}_n \right]
      -V\left[ \epsilon^*_{n+1}\phi^*_n+\epsilon^*_n\phi^*_{n+1} \right] ,
\end{align}
and the expansions of the polynomials and associated functions as given in 
(\ref{phiexp:a}-\ref{epsexp:d}).
\end{proof}

In addition to the coupled equations of Proposition \ref{Linear1} and Corollary \ref{Linear2}, 
evaluations of the spectral coefficients at the singular points satisfy bilinear relations.
\begin{proposition}[\cite{FW_2006a}]
For $ j=1,\ldots,M $ (i.e. $ z_j \neq 0 $) the evaluations of the spectral coefficients satisfy the 
recurrence and functional relations  
\begin{gather}
  \Omega^2_n(z_j)
  = \frac{\kappa_n \phi_{n+2}(0)}{\kappa_{n+1} \phi_{n+1}(0)}z_j
      \Theta_n(z_j)\Theta_{n+1}(z_j)+V^2(z_j) ,
  \label{OTeq:a} \\
  \Omega^{*2}_n(z_j)
  = \frac{\kappa_n \bar{\phi}_{n+2}(0)}{\kappa_{n+1} \bar{\phi}_{n+1}(0)}z_j
      \Theta^*_n(z_j)\Theta^*_{n+1}(z_j)+V^2(z_j) ,
  \label{OTeq:b} \\
  \left[\Omega_{n-1}(z_j)
    -\frac{\kappa^2_{n-1}}{\kappa^2_{n}}\frac{\phi_{n+1}(0)}{\phi_{n}(0)}\Theta_n(z_j)
  \right]^2
  = \frac{\phi_{n+1}(0)\bar{\phi}_{n}(0)}{\kappa^2_{n}}
      \Theta_n(z_j)\Theta^*_{n-1}(z_j)+V^2(z_j) ,
  \label{OTeq:c}
\end{gather}
\begin{equation}
  \left[\Omega^*_{n-1}(z_j)
    -\frac{\kappa^2_{n-1}}{\kappa^2_{n}}\frac{\bar{\phi}_{n+1}(0)}{\bar{\phi}_{n}(0)}
      z_j\Theta^*_n(z_j)
  \right]^2
  = \frac{\kappa_{n-1}\bar{\phi}_{n+1}(0)\phi_{n}(0)}{\kappa^3_{n}}
      z^2_j\Theta^*_n(z_j)\Theta_{n-1}(z_j)+V^2(z_j) ,
  \label{OTeq:d}
\end{equation}
\begin{equation}
    \frac{\phi_{n+1}(0)\bar{\phi}_{n+1}(0)}{\kappa^2_{n}}
     z_j\Theta_n(z_j)\Theta^*_{n}(z_j)
  = \left[\Omega_{n}(z_j)+V(z_j)-\frac{\kappa_{n+1}}{\kappa_{n}}z_j\Theta_n(z_j)\right]
    \left[\Omega^*_{n}(z_j)-V(z_j)-\frac{\kappa_{n+1}}{\kappa_{n}}\Theta^*_n(z_j)\right] .
\label{OTeq:e}
\end{equation}
\end{proposition}
It is these relations that lead directly to one of the pair of coupled discrete 
Painlev\'e equations. 

We note the initial members of the sequences of spectral coefficients
$ \{\Theta_n\}^{\infty}_{n=0} $, $ \{\Theta^*_n\}^{\infty}_{n=0} $,
$ \{\Omega_n\}^{\infty}_{n=0} $, $ \{\Omega^*_n\}^{\infty}_{n=0} $ are given by
\begin{align}
  2\frac{\phi_1(0)}{\kappa_0}\Theta_0(z) = & 2V(z)-\kappa^2_0 U(z) ,
  \label{Theta0} \\
  2\frac{\bar{\phi}_1(0)}{\kappa_0}z\Theta^*_0(z) = & -2V(z)-\kappa^2_0 U(z) ,
  \label{ThetaS0} \\ 
  2\phi_1(0)\Omega_0(z) = & \kappa_1z[2V(z)-\kappa^2_0 U(z)] - \kappa^2_0\phi_1(0)U(z) ,
  \label{Omega0} \\
  2\bar{\phi}_1(0)z\Omega^*_0(z) = & -\kappa_1[2V(z)+\kappa^2_0 U(z)] - \kappa^2_0\bar{\phi}_1(0)zU(z) ,
  \label{OmegaS0}
\end{align}
and observe that $ U(z) $ defines the initial values for the recurrences of 
Proposition \ref{Linear1}.

Hereafter we restore a singularity of the weight at $ z=0 $ (i.e. $ \rho_0 \neq 0 $) 
in addition to the previous finite ones so that the degrees of $ W, V $ are augmented by a unit
(now $ M+1 $ and $ M $ respectively) and the spectral coefficients $ \Theta_n, \Theta^*_n $ 
and $ \Omega_n, \Omega^*_n $ are polynomials of degree $ M-1 $ and $ M $ respectively.
In this setting the spectral matrix $ A_n $ has a partial fraction decomposition
\begin{equation}
   A_n = \sum^{M}_{j=0}\frac{A_{n,j}}{z-z_j}, \quad z_0=0,
\label{AnPF}
\end{equation}
which define the residue matrices $ A_{n,j} $for $ j=0,\ldots, M $. We remark that the spectral
matrix has singularities at $ z=0 $ and $ z=\infty $ irregardless of the locations of singularities
of the weight (i.e. zeros of $ W $) due to the bi-orthogonality structure.
The $j$-th residue matrix $ A_{n,j} $ at the finite singularity $ z_j $ is given by
\begin{align}
   A_{n,j} & = \frac{1}{W'(z_j)}
       \begin{pmatrix}
              -\Omega_n(z_j)-V(z_j)
              +\dfrac{\kappa_{n+1}}{\kappa_n}z_j\Theta_n(z_j)
            & \dfrac{\phi_{n+1}(0)}{\kappa_n}\Theta_n(z_j)
            \cr
              -\dfrac{\bar{\phi}_{n+1}(0)}{\kappa_n}z_j\Theta^*_n(z_j)
            &  \Omega^*_n(z_j)-V(z_j)
                     -\dfrac{\kappa_{n+1}}{\kappa_n}\Theta^*_n(z_j)
            \cr
       \end{pmatrix} ,
\label{An_resJ}
\end{align}
with $ 2V(z_j)=\rho_j W'(z_j) $ for $ j=1,\ldots,M $, while for $ j=0 $ the expression is
\begin{equation}
   A_{n,0} = (n-\rho_0)
       \begin{pmatrix}
              1 & -r_n \cr 0 & 0 \cr
       \end{pmatrix} ,
\label{An_res0}
\end{equation}
and for the singular point $ z_{M+1}=\infty $ it is
\begin{equation}
   A_{n,\infty} = -\sum^{M}_{j=0}A_{n,j} =
       \begin{pmatrix}
              -n & 0 \cr -(n+\sum^{M}_{j=0}\rho_j)\bar{r}_n & \sum^{M}_{j=0}\rho_j \cr
       \end{pmatrix} .
\label{An_resInfty}
\end{equation}
In the following section we will draw heavily on partial fraction decompositions of the
spectral coefficients which imply the summation identities
\begin{align}
  \sum^M_{j=0} \frac{\Theta_n(z_j)}{W'(z_j)}
 & = 0 ,
\\
  \sum^M_{j=0} \frac{\Theta^*_n(z_j)}{W'(z_j)}
 & = -(n+\sum^M_{j=0}\rho_j)\frac{\bar{\phi}_n(0)}{\bar{\phi}_{n+1}(0)} ,
\\
  \sum^M_{j=0} \frac{\Omega_n(z_j)-V(z_j)-\frac{\DySt \kappa_{n+1}}{\DySt \kappa_n}z_j\Theta_n(z_j)}{W'(z_j)}
 & = -(n+\sum^M_{j=0}\rho_j) ,
\\
  \sum^M_{j=0} \frac{\Omega^*_n(z_j)-V(z_j)-\frac{\DySt \kappa_{n+1}}{\DySt \kappa_n}\Theta^*_n(z_j)}{W'(z_j)}
 & = -\sum^M_{j=0}\rho_j .
\end{align}

For simplicity we parameterise the free singularities $ z_j(t) $ as arbitrary trajectories with
respect to a single deformation variable $ t $ so that 
\begin{equation}
  \frac{d}{dt} = \sum^M_{j=0} \dot{z}_j\frac{\partial}{\partial z_j} ,
\end{equation}
where we include $ j=0 $ in the sum for convenience even though $ \dot{z}_0=0 $.
\begin{proposition}[\cite{FW_2006a}]\label{BOPS_deformDer}
The deformation derivatives of the system (\ref{Ydefn}) with respect to arbitrary 
deformations of the singularities $ z_j $ are given by
\begin{equation}
   \frac{d}{dt}Y_{n} = B_{n} Y_{n}
   = \left[ B_{\infty} - \sum^M_{j=0}\dot{z}_j\frac{A_{n,j}}{z-z_j} \right] Y_{n} ,
\label{YdeformDer}
\end{equation}
where
\begin{equation}
   B_{\infty} =   \begin{pmatrix}
    \frac{\DySt\dot{\kappa}_n}{\DySt\kappa_n} & 0 \cr
    \frac{\DySt\dot{\kappa}_n\bar{\phi}_n(0)+\kappa_n\dot{\bar{\phi}}_n(0)}{\DySt\kappa_n} & -\frac{\DySt\dot{\kappa}_n}{\DySt\kappa_n} \cr
  \end{pmatrix} .
\end{equation}
\end{proposition}
A particularly important deduction from (\ref{YdeformDer}) are the dynamics of the $r$-coefficients
\begin{align}
   \frac{\dot{r}_n }{r_n} & =
   \sum^M_{j=0}\dot{z}_j\frac{\Omega_{n-1}(z_j)-V(z_j)}{W'(z_j)} ,
   \label{rdot} \\
   \frac{\dot{\bar{r}}_n}{\bar{r}_n} & =
   \sum^M_{j=0}\dot{z}_j\frac{\Omega^*_{n-1}(z_j)+V(z_j)}{W'(z_j)} .
   \label{rCdot}
\end{align}

Compatibility of the spectral derivative (\ref{YspecDer}) and the deformation 
derivative (\ref{YdeformDer}) leads to the Schlesinger equations for the residue matrices.
\begin{proposition}[\cite{FW_2006a}]
The residue matrices satisfy a system of integrable, non-linear partial differential 
equations, the Schlesinger equations, 
\begin{gather}
   \dot{A}_{n,j} = \left[ B_{\infty},A_{n,j} \right]
   + \sum_{{\ScSt k\neq j}\atop{\ScSt 0\leq k \leq M}} \frac{\dot{z}_j-\dot{z}_k}{z_j-z_k} \left[ A_{n,k}, A_{n,j} \right],
   \qquad j=0,\ldots, M , 
   \label{SchlesingerEqn}\\
   \dot{A}_{n,\infty} = \left[ B_{\infty},A_{n,\infty} \right] .
\label{An_SE}
\end{gather} 
\end{proposition}

In the following section we will require a more explicit form for the Schlesinger equations,
which were first given in \cite{FW_2006a}. 
\begin{lemma}
The deformation derivatives of the residues of the spectral matrix are given in 
component form by
\begin{multline}
   \frac{\kappa_n}{\phi_{n+1}(0)}W'(z_j)\dot{A}_{nj,12}
   = \frac{2\dot{\kappa}_n}{\kappa_n}\Theta_n(z_j) 
  \\
   +\sum_{{\ScSt k\neq j}\atop{\ScSt 0\leq k \leq M}}
     \frac{1}{W'(z_k)}\frac{\dot{z}_j-\dot{z}_k}{z_j-z_k}
     \left\{ 
       \Theta_n(z_k)\left[ 2\Omega_n(z_j)-2\frac{\kappa_{n+1}}{\kappa_n}z_j\Theta_n(z_j)+n\frac{W(z_j)}{z_j} \right]
     \right.
  \\ \left.
      -\Theta_n(z_j)\left[ 2\Omega_n(z_k)-2\frac{\kappa_{n+1}}{\kappa_n}z_k\Theta_n(z_k)+n\frac{W(z_k)}{z_k} \right]
     \right\} ,
\label{AnSE:a}
\end{multline}
and
\begin{multline}
   W'(z_j)\dot{A}_{nj,11}
   = -\frac{\phi_{n+1}(0)}{\kappa_n}\frac{\dot{\kappa}_n\bar{\phi}_n(0)+\kappa_n\dot{\bar{\phi}}_n(0)}{\kappa^2_n}\Theta_n(z_j) 
  \\
   +\sum_{{\ScSt k\neq j}\atop{\ScSt 0\leq k \leq M}}
     \frac{1}{W'(z_k)}\frac{\dot{z}_j-\dot{z}_k}{z_j-z_k}
     \left\{ 
       \frac{\Theta_n(z_j)}{\Theta_n(z_k)}
  \left[ \Omega_n(z_k)+V(z_k)-\frac{\kappa_{n+1}}{\kappa_n}z_k\Theta_n(z_k) \right]
  \left[ \Omega_n(z_k)-V(z_k)-\frac{\kappa_{n+1}}{\kappa_n}z_k\Theta_n(z_k)+n\frac{W(z_k)}{z_k} \right]
     \right.
  \\ \left.
      -\frac{\Theta_n(z_k)}{\Theta_n(z_j)}
  \left[ \Omega_n(z_j)+V(z_j)-\frac{\kappa_{n+1}}{\kappa_n}z_j\Theta_n(z_j) \right]
  \left[ \Omega_n(z_j)-V(z_j)-\frac{\kappa_{n+1}}{\kappa_n}z_j\Theta_n(z_j)+n\frac{W(z_j)}{z_j} \right]
     \right\} ,
\label{AnSE:b}
\end{multline}
for $ j=0,\ldots,M $.
\end{lemma}
\begin{proof}
These follow from the Schlesinger equations (\ref{SchlesingerEqn}) and the transition formulae
\begin{gather}
  \Omega^*_n(z_j)-\frac{\kappa_{n+1}}{\kappa_n}\Theta^*_n(z_j)
 = \Omega_n(z_j)-\frac{\kappa_{n+1}}{\kappa_n}z_j\Theta_n(z_j)+n\frac{W(z_j)}{z_j} ,
\label{Tform:a}
\\
  \frac{\phi_{n+1}(0)\bar{\phi}_{n+1}(0)}{\kappa^2_n}z_j\Theta^*_n(z_j)
  =   \frac{1}{\Theta_n(z_j)}
  \left[ \Omega_n(z_j)+V(z_j)-\frac{\kappa_{n+1}}{\kappa_n}z_j\Theta_n(z_j) \right]
  \left[ \Omega_n(z_j)-V(z_j)-\frac{\kappa_{n+1}}{\kappa_n}z_j\Theta_n(z_j)+n\frac{W(z_j)}{z_j} \right] ,
\label{Tform:b}
\end{gather}
for $ j=0,\ldots,M $ and where the ratio $ W(z_j)/z_j $ is interpreted as the limit of $ z_0 \to 0 $ 
for $ j=0 $. 
\end{proof}

For each singular point $ z_j $ the monodromy matrix $ M_j $ is defined by 
\begin{equation}
   Y_n(z_j + \delta e^{2\pi i}) = Y_n(z_j + \delta)M_j ,
\end{equation}
and has the classical, triangular structure
\begin{equation}
   M_j =
       \begin{pmatrix}
              1 & c_j(1-e^{-2\pi i\rho_j}) \cr
              0 & e^{-2\pi i\rho_j} \cr
       \end{pmatrix}, \quad j=0,\ldots,M ,
\end{equation}
where $ c_j $ is independent of the $ z_j $, and thus $ t $, and $ n \in \Z_{\geq 0} $ 
but depends on other details of the weight.

\section{The Hamiltonian Formulation and Garnier Systems}\label{Hamiltonian}
\setcounter{equation}{0}

We fix the singularities in canonical position (see \cite{IKSY_1991}), i.e. taking them at distinct
points
\begin{equation}
   z_0 = 0, z_1, \ldots, z_N, z_{N+1}=1, z_{N+2}=\infty, 
\end{equation}
with exponents $ \rho_0, \rho_1, \ldots, \rho_N, \rho_{N+1}=\rho, \rho_{N+2}=\rho_{\infty} $ 
respectively, so that the number of finite singularities is $ N+2 $. Note that now we have a
singularity of the weight at the origin, i.e. $ \rho_0 \neq 0 $.
The denominator polynomial for the weight data is given by 
\begin{equation}
   W(z) = z(z-1)\prod^N_{j=1}(z-z_j)
        = z\sum^{N+1}_{l=0} (-)^{N+1-l}z^l e_{N+1-l} ,
\label{Wdefn}
\end{equation}
where the elementary symmetric functions of the singularity positions are denoted 
$ e_l, l=0, \ldots,N+1 $ and in particular $ e_0=1 $, $ e_{N+2}=0 $ and $ e_{N+1}=\prod^N_{j=1}z_j $.
The numerator polynomial is 
\begin{equation}
   2V(z) = z(z-1)\prod^N_{j=1}(z-z_j)
             \left\{ \frac{\rho_0}{z}+\frac{\rho}{z-1}+\sum^N_{j=1}\frac{\rho_j}{z-z_j} \right\}
         = \sum^{N+1}_{l=0} (-)^{N+1-l}z^l m_{N+1-l} ,
\label{2Vdefn}
\end{equation}
where the last relation defines the coefficients $ m_l, l=0,\ldots,N+1 $ and 
we observe that $ m_0=\rho_0+\rho+\sum^N_{j=1}\rho_j $ and $ m_{N+1}=\rho_0 e_{N+1} $.
We can parameterise the upper off-diagonal element of the spectral matrix, i.e. the spectral
coefficient $ \Theta_n $, which is now of degree $ N $, so that
\begin{equation}
 \Theta_n(z) = \Theta_{\infty}\prod^{N}_{r=1}(z-q_r) , \quad
 \Theta_{\infty} = (n+1+m_0)\frac{\kappa_n}{\kappa_{n+1}} ,
\label{ThetaRep}
\end{equation}
and thus
\begin{equation}
	\frac{\Theta'_n}{\Theta_n} = \sum^{N}_{r=1} \frac{1}{z-q_r} ,
\end{equation}
where the poles $ q_r $ will play the role of canonical co-ordinates and analogue of the sixth 
Painlev\'e transcendent.
For notational simplicity we will often suppress the $ n $ index dependency as we 
do not discuss recurrences in this variable in this section.
Furthermore we will assume throughout {\it generic conditions} on the dependent and independent variables,
namely that
\begin{enumerate}
  \item[(i)]
   non-coincidence of singular points, $ z_j \neq z_k $ for $ j\neq k $ and $ j,k=0,\ldots, N+2 $,
   so that $ W'(z_j) \neq 0 $ for $ j=0,\ldots, N+2 $
  \item[(ii)]
   avoidance by the co-ordinates with each other $ q_r \neq q_s $ for $ r\neq s $ and {$ r,s=1,\ldots, N $}
   implying that $ \Theta'_n(q_r) \neq 0 $ for $ r=1,\ldots, N $
   and with the fixed singularities $ q_r \neq z_j $ for $ r=1\ldots,N $ and $ j=0,\ldots,N+2 $
   and consequently $ W(q_r) \neq 0 $ for $ r=1,\ldots, N $ and $ \Theta_n(z_j) \neq 0 $ for $ j=0,\ldots,N+2 $.
\end{enumerate}

Prior to stating our main result we note some summation identities for $ j=1,\ldots,N $
\begin{gather}
   \sum_{{\ScSt 0\leq k \leq N+1}\atop{\ScSt k\neq j}}\frac{1}{z_j-z_k} 
  = \frac{1}{2}\frac{W''(z_j)}{W'(z_j)} ,
\label{Ssum:a} \\
   \sum_{{\ScSt 0\leq k \leq N+1}\atop{\ScSt k\neq j}}\frac{\rho_k}{z_j-z_k} 
  = \frac{2V'(z_j)}{W'(z_j)}-\frac{V(z_j)W''(z_j)}{[W'(z_j)]^2} ,   
\label{Ssum:b} \\
   \sum_{{\ScSt 0\leq k \leq N+1}\atop{\ScSt k\neq j}}\frac{\Theta_n(z_k)}{W'(z_k)}\frac{1}{z_j-z_k} 
  = \frac{\Theta'_n(z_j)}{W'(z_j)}-\frac{1}{2}\frac{\Theta_n(z_j)W''(z_j)}{[W'(z_j)]^2} ,
\label{Ssum:c}
\end{gather}
and
\begin{equation}
   \sum_{{\ScSt 1\leq s \leq N}\atop{\ScSt s\neq r}}\frac{1}{q_r-q_s} 
  = \frac{1}{2}\frac{\Theta''_n(q_r)}{\Theta'_n(q_r)} . 
\label{Tsum:a}
\end{equation}
In addition partial fraction expansions imply the following summation identities
\begin{equation}
  \sum^{N+1}_{j=0} \frac{\Theta_n(z_j)}{W'(z_j)}z^{\sigma}_j =
  \begin{cases}
     0, & \sigma = 0 \\
     \Theta_{\infty}, & \sigma = 1 \\
     \Theta_{\infty}[1+\sum^N_{k=1}z_k-\sum^N_{r=1}q_r], & \sigma = 2
  \end{cases} ,
\label{Ssum:d}
\end{equation}
\begin{equation}
  \sum^{N+1}_{j=0} \frac{\Theta_n(z_j)}{W'(z_j)}\frac{z^{\sigma}_j}{z_j-q_r} =
  \begin{cases}
     0, & \sigma = 0,1 \\
     \Theta_{\infty}, & \sigma = 2 \\
     \Theta_{\infty}[1+\sum^N_{k=1}z_k-\sum^N_{s\neq r}q_s], & \sigma = 3
  \end{cases} ,
\label{Ssum:e}
\end{equation}
for $ r=1,\ldots, N $,
\begin{equation}
  \sum^{N+1}_{j=0} \frac{\Theta_n(z_j)}{W'(z_j)}\frac{z^{\sigma}_j}{(z_j-q_r)(z_j-q_s)}
  = \begin{cases}
     -\delta_{s,r}q^{\sigma}_r\frac{\DySt \Theta'_n(q_r)}{\DySt W(q_r)}, & \sigma = 0,1,2 \\
     \Theta_{\infty}-\delta_{s,r}q^{3}_r\frac{\DySt \Theta'_n(q_r)}{\DySt W(q_r)}, & \sigma = 3 \\
     \Theta_{\infty}[1+\sum^N_{k=1}z_k-\sum^N_{t=1}q_t+q_r+q_s]-\delta_{s,r}q^{4}_r\frac{\DySt \Theta'_n(q_r)}{\DySt W(q_r)}, & \sigma = 4
    \end{cases} ,
\label{Ssum:f}
\end{equation}
for $ r,s=1,\ldots, N $, and
\begin{multline}
  \sum^{N+1}_{j=0} \frac{\Theta_n(z_j)}{W'(z_j)}\frac{z^{\sigma}_j}{(z_j-q_r)(z_j-q_s)(z_j-q_t)}
  \\
  = -(1-\delta_{t,s})\delta_{r,s}q^{\sigma}_s\frac{\Theta'_n(q_s)}{(q_r-q_t)W(q_s)}
    -(1-\delta_{s,r})\delta_{t,r}q^{\sigma}_r\frac{\Theta'_n(q_r)}{(q_t-q_s)W(q_r)}
    -(1-\delta_{r,t})\delta_{s,t}q^{\sigma}_t\frac{\Theta'_n(q_t)}{(q_s-q_r)W(q_t)}
  \\
    +\delta_{r,s}\delta_{s,t}q^{\sigma}_r\frac{\Theta'_n(q_r)}{W(q_r)}
     \left[ \frac{W'(q_r)}{W(q_r)}-\frac{1}{2}\frac{\Theta''_n(q_r)}{\Theta'_n(q_r)}-\frac{\sigma}{q_r} \right] ,
\label{Ssum:g}
\end{multline}
for $ r,s,t=1,\ldots, N $ with $ \sigma=0,1,2,3 $ whilst for $ \sigma=4 $ an additional 
term $ \Theta_{\infty} $ is necessary.

\begin{proposition}
Assume that generic conditions apply.
The dynamics of the bi-orthogonal system is governed by the Hamiltonian dynamics of 
the Garnier system $ {\mathcal G}_N \equiv \{q_j,p_j;K_j,z_j\} $ with co-ordinate 
$ q_r $ defined above and momenta $ p_r $ given by
\begin{equation}
     p_r = -\frac{\Omega_n(q_r)+V(q_r)}{W(q_r)}=A_{1,1}(q_r), \quad r=1,\ldots, N ,
\label{Hamp}
\end{equation}
and the Hamiltonian by
\begin{equation}
  K_j = \frac{\Theta_n(z_j)}{W'(z_j)}\sum^N_{r=1}\frac{W(q_r)}{\Theta'_n(q_r)}\frac{1}{z_j-q_r}
    \left[ p^2_r+p_r\left( \frac{2V(q_r)}{W(q_r)}-\frac{n}{q_r}-\frac{1}{z_j-q_r} \right)
      -\frac{n(1+m_0)}{q_r(q_r-1)} 
    \right] .
\label{Ham}
\end{equation}
The indicial exponents are given by $ \theta_j = -\rho_j $ for $ j=1,\ldots, N+1 $,
$ \theta_0 = n-\rho_0 $, and $ \alpha_{\infty} = -n, \theta_{\infty} = n+1+\sum^{N+1}_{j=0}\rho_j $ 
with the constant $ \kappa = -n(1+m_0) $.  The latter relation is the
necessary condition for a classical solution to the Garnier system, and for the "seed" 
solution of $ n=0 $ it vanishes.
\end{proposition}
\begin{proof}
We seek to compare our system with the second order ODE for the Garnier system Eqs. (4.1.1) 
and (4.1.4) in \cite{IKSY_1991}, see also \cite{Ga_1912}, \cite{KO_1984},
\begin{equation}
	\phi'' 
	+\left\{ \sum^{N+1}_{j=0}\frac{1-\theta_j}{z-z_j}-\sum^{N}_{r=1}\frac{1}{z-q_r} \right\}\phi'
	+\left\{\frac{\kappa}{z(z-1)}-\sum^{N}_{j=1}\frac{z_j(z_j-1)K_j}{z(z-1)(z-z_j)}
	               +\sum^{N}_{r=1}\frac{q_r(q_r-1)p_r}{z(z-1)(z-q_r)}
	 \right\}\phi = 0 .
\label{garnierODE}
\end{equation}
Using Prop. \ref{scalarODE} we find that
\begin{equation}
	p_1 = \frac{\rho_0+1-n}{z}+\frac{\rho+1}{z-1}
	      +\sum^{N}_{j=1}\frac{\rho_j+1}{z-z_j}
	      -\sum^{N}_{j=1}\frac{1}{z-q_j} .
\end{equation}
From the residues of $ p_1 $ at $ z=z_j $ for $ j = 0,1,\ldots, N+1 $ we can read off 
the indicial exponents. We also note that 
\begin{equation}
	p_2 = \frac{1}{z(z-1)}\left[ -n(1+m_0)+{\rm O}(z^{-1}) \right] ,
\end{equation}
as $ z \to \infty $.

An alternative form of $ p_2 $ to (\ref{2ODEcff:b}) can be given in terms of the 
residue matrices 
\begin{equation}
   p_2 = \sum^{N+1}_{j=0}\sum^N_{r=1}\frac{A_{nj,11}}{(z-z_j)(z-q_r)}
  \\
  + \sum_{0\leq j<k\leq N+1}\frac{{\rm Tr}A_{nj}{\rm Tr}A_{nk}-{\rm Tr}A_{nj}A_{nk}-A_{nj,11}-A_{nk,11}}{(z-z_j)(z-z_k)} .
\end{equation}
From $ p_r = {\rm Res}_{z=q_r}p_2(z) $, the above relation and (\ref{An_resJ}) we 
find (\ref{Hamp}). We need to invert this relationship and express $ \Omega_n+V $ in 
terms of the canonical co-ordinate and momenta. Because this variable is a polynomial 
of degree $ N+1 $ in the spectral variable $ z $ we can use the Lagrange interpolation 
formula at the nodes $ z=0,q_1,\ldots,q_N, 1 $ (which are assumed to be distinct)
\begin{multline}
  \Omega_n+V = \sum^N_{r=1} \left[ \Omega_n(q_r)+V(q_r) \right]
                            \frac{z(z-1)\prod_{s\neq r}(z-q_s)}{q_r(q_r-1)\prod_{s\neq r}(q_r-q_s)}
  \\
      +  \left[ \Omega_n(0)+V(0) \right]
                            \frac{(1-z)\prod_{s}(z-q_s)}{\prod_{s}(-q_s)}
      + \left[ \Omega_n(1)+V(1) \right]
                            \frac{z\prod_{s}(z-q_s)}{\prod_{s}(1-q_s)} .
\end{multline}
The coefficients of the two last terms are known as 
$ \Omega_n(0)+V(0) = (-)^N(n-\rho_0)e_{N+1} $ from (\ref{Omexp:b}) and utilising 
the coefficient of $ z^{N+1} $ in (\ref{Omexp:a}) we deduce 
\begin{equation}
  \frac{\Omega_n(1)+V(1)}{\Theta_n(1)} = \frac{1+m_0}{\Theta_{\infty}}+(-)^N(n-\rho_0)\frac{e_{N+1}}{\Theta_N(0)}
 + \sum^N_{r=1}\frac{p_rW(q_r)}{q_r(q_r-1)\Theta'_n(q_r)} .
\end{equation}
Consequently we conclude that
\begin{equation}
  \Omega_n(z)+V(z)-\frac{\kappa_{n+1}}{\kappa_n}z\Theta_n(z)
   = \Theta_n(z)\left[ -\frac{n}{\Theta_{\infty}}z +(-)^N(n-\rho_0)\frac{e_{N+1}}{\Theta_n(0)}
         - \sum^N_{r=1}\frac{z}{z-q_r}\frac{p_rW(q_r)}{q_r\Theta'_n(q_r)} \right] .
\label{OmegaRep}
\end{equation}

We also require a similar representation of $ 2V(z) $ and proceeding in the same 
manner we find
\begin{equation}
  2V(z) = \Theta_n(z)
  \left[ -\rho_0\frac{W'(0)}{\Theta_{n}(0)}(z-1) +\rho\frac{W'(1)}{\Theta_n(1)}z
          +z(z-1)\sum^N_{r=1}\frac{1}{z-q_r}\frac{2V(q_r)}{q_r(q_r-1)\Theta'_n(q_r)} \right] .
\label{2VRep}
\end{equation}
From this, under the limit $ z \to \infty $, we have the summation
\begin{equation}
  \sum^N_{r=1}\frac{2V(q_r)}{q_r(q_r-1)\Theta'_n(q_r)} 
 = \frac{m_0}{\Theta_{\infty}}+\rho_0\frac{W'(0)}{\Theta_{n}(0)}-\rho\frac{W'(1)}{\Theta_n(1)} ,
\label{Tsum:b}
\end{equation}
and for $ z=z_j $ we have
\begin{equation}
   z_j(z_j-1)\sum^N_{r=1}\frac{2V(q_r)}{q_r(q_r-1)\Theta'_n(q_r)}\frac{1}{z_j-q_r} 
   = \rho_0\frac{W'(0)}{\Theta_n(0)}(z_j-1)-\rho\frac{W'(1)}{\Theta_{n}(1)}z_j+\frac{2V(z_j)}{\Theta_n(z_j)} .
\label{Tsum:h}
\end{equation}
A number of other sums for $ j=1,\ldots,N $ which will be required subsequently
are
\begin{equation}
  \sum^N_{r=1}\frac{2V(q_r)}{(z_j-q_r)^2\Theta'_n(q_r)} 
 = \frac{m_0}{\Theta_{\infty}}-\frac{2V'(z_j)}{\Theta_{n}(z_j)}+\frac{2V(z_j)\Theta'_n(z_j)}{\Theta^2_n(z_j)} ,
\label{Tsum:c}
\end{equation}
\begin{equation}
  \sum^N_{r=1}\frac{2V(q_r)}{(z_j-q_r)q_r\Theta'_n(q_r)} 
 = -\frac{m_0}{\Theta_{\infty}}-\frac{2V'(0)}{z_j\Theta_{n}(0)}+\frac{2V(z_j)}{z_j\Theta_n(z_j)} ,
\label{Tsum:d}
\end{equation}
and, for $ r=1,\ldots,N $ the sum
\begin{multline}
   q_r(q_r-1)\sum^N_{s\neq r}\frac{2V(q_s)}{q_s(q_s-1)\Theta'_n(q_s)}\frac{1}{q_r-q_s} 
  \\
   = \rho_0\frac{W'(0)}{\Theta_n(0)}(q_r-1)-\rho\frac{W'(1)}{\Theta_{n}(1)}q_r
  +\frac{2V(q_r)}{\Theta'_n(q_r)}\left[
         \frac{V'(q_r)}{V(q_r)}-\frac{1}{2}\frac{\Theta''_n(q_r)}{\Theta'_n(q_r)}-\frac{2q_r-1}{q_r(q_r-1)} \right] .
\label{Tsum:e}
\end{multline}

The analogous representation for $ W(z) $ takes the form
\begin{equation}
  W(z) = z(z-1)\Theta_n(z)
  \left[ \frac{1}{\Theta_{\infty}}+\sum^N_{r=1}\frac{1}{z-q_r}\frac{W(q_r)}{q_r(q_r-1)\Theta'_n(q_r)} \right] ,
\label{WRep}
\end{equation}
which enables one to infer the summations
\begin{equation}
  \sum^N_{r=1}\frac{W(q_r)}{(z_j-q_r)q_r(q_r-1)\Theta'_n(q_r)} = -\frac{1}{\Theta_{\infty}} ,
  \quad j=1,\ldots,N ,
\label{Tsum:f}
\end{equation}
and for $ r=1,\ldots,N $ the sum
\begin{equation}
  \sum^N_{s\neq r}\frac{W(q_s)}{q_s(q_s-1)\Theta'_n(q_s)}\frac{1}{q_r-q_s} 
   = -\frac{1}{\Theta_{\infty}}
  +\frac{W(q_r)}{q_r(q_r-1)\Theta'_n(q_r)}\left[
         \frac{W'(q_r)}{W(q_r)}-\frac{1}{2}\frac{\Theta''_n(q_r)}{\Theta'_n(q_r)}-\frac{2q_r-1}{q_r(q_r-1)} \right] .
\label{Tsum:g}
\end{equation}

We proceed in two steps - the first being to compute the derivatives of the 
canonical variables using the theory from \cite{FW_2006a} and the second to compute 
the Hamiltonian itself. Then all we require is a verification of the Hamilton equations 
of motion.

Our first task is the computation of the deformation derivative of $ q_r $
\begin{equation}
   \dot{q}_r = \sum^N_{k=1} \dot{z}_k \frac{\partial q}{\partial z_k}
  = -{\rm Res}_{z=q_r}\frac{\dot{A}_{12}}{A_{12}} ,
\end{equation}
and where we use the partial fraction expansion (\ref{AnPF}) for the $ 1,2 $ element
and the Schlesinger equation (\ref{AnSE:a}) for $ \dot{A}_{nj,12} $. Now we employ 
the expressions for $ \Theta_n(z) $ and 
$ \Omega_n(z)-\frac{\DySt \kappa_{n+1}}{\DySt \kappa_n}z\Theta_n(z) $ at $ z=z_j,z_k $ 
in terms of the canonical variables as given by (\ref{ThetaRep}) and (\ref{OmegaRep}) 
along with (\ref{2VRep}) in this formula. After interchanging the summation order 
each term contains a factor which is a sum over the singularity locations and we can 
evaluate these using the $ \sigma=0 $ cases of (\ref{Ssum:d}), (\ref{Ssum:e}) and 
(\ref{Ssum:f}). Simplifying we find for $ r,j=1,\ldots,N $
\begin{equation}
  (z_j-q_r)\frac{\partial}{\partial z_j}q_r 
  = \frac{\Theta_n(z_j)W(q_r)}{\Theta'_n(q_r)W'(z_j)}
    \left[ 2p_r+\frac{2V(q_r)}{W(q_r)}-\frac{n}{q_r}-\frac{1}{z_j-q_r} \right] .
\label{Ham_qDer}
\end{equation}

We now proceed to the computation of $ \dot{p}_r $ which is a more laborious task. 
We start with the partial fraction expansion for $ A_{n,11}(q_r) $ given by (\ref{AnPF}), 
differentiate with respect to $ t $ and employ the $ 1,1 $ component of the Schlesinger 
equation (\ref{AnSE:b}) for $ \dot{A}_{nj,11} $. Into this expression we must substitute 
the expressions for $ \Theta_n(z) $ and 
$ \Omega_n(z)\pm V(z)-\frac{\DySt \kappa_{n+1}}{\DySt \kappa_n}z\Theta_n(z) $ at 
$ z=z_j,z_k $ in terms of the canonical variables as given by (\ref{ThetaRep}) and 
(\ref{OmegaRep}) along with (\ref{2VRep}). Again we interchange the summation order 
which yields in each term a factor of a singularity sum. To evaluate these sums we 
employ all of the formulae given in (\ref{Ssum:d}), (\ref{Ssum:e}), (\ref{Ssum:f}) 
and (\ref{Ssum:g}) and assemble all the terms together. We also require the formula
for $ \dot{q}_r $ implied by (\ref{Ham_qDer}). Considerable cancellation occurs at 
this stage but we are not yet finished. To arrive at the final result we must employ 
the transcendent sums (\ref{Tsum:b}) and (\ref{Tsum:e}), and the result for 
$ r,j=1,\ldots,N $ is 
\begin{multline}
  (z_j-q_r)\frac{W'(z_j)}{\Theta_n(z_j)}\frac{\partial}{\partial z_j}p_r
 = -\frac{W(q_r)}{\Theta'_n(q_r)}
   \left[
    p^2_r\left( \frac{W'(q_r)}{W(q_r)}-\frac{1}{2}\frac{\Theta''_n(q_r)}{\Theta'_n(q_r)} \right)
   +p_r\frac{2V(q_r)}{W(q_r)}\left( \frac{V'(q_r)}{V(q_r)}-\frac{1}{2}\frac{\Theta''_n(q_r)}{\Theta'_n(q_r)} \right)
   \right.
 \\ \left.
   -n\frac{p_r}{q_r}\left( \frac{W'(q_r)}{W(q_r)}-\frac{1}{2}\frac{\Theta''_n(q_r)}{\Theta'_n(q_r)}-\frac{1}{q_r} \right)
   -\frac{p_r}{z_j-q_r}\left( \frac{W'(q_r)}{W(q_r)}-\frac{1}{2}\frac{\Theta''_n(q_r)}{\Theta'_n(q_r)}+\frac{1}{z_j-q_r} \right)
   +n(1+m_0)\frac{1}{q_r(q_r-1)(z_j-q_r)}
   \right]                     
 \\
  -\sum_{s\neq r}\frac{W(q_s)}{\Theta'_n(q_s)}\frac{1}{q_s-q_r}
   \left[
     p^2_s+p_s\frac{2V(q_s)}{W(q_s)}-n\frac{p_s}{q_s}
    -\frac{p_s}{z_j-q_s}+n(1+m_0)\frac{q_s-q_r}{q_s(q_s-1)(z_j-q_s)}
   \right] .
\label{Ham_pDer}
\end{multline}

Our starting point for the computation of the Hamiltonian is the formula
\begin{equation}
   K_j = -A_{nj,11}\sum^N_{r=1}\frac{1}{z-q_r}
  - \sum_{{\ScSt 0\leq k\leq N+1}\atop{\ScSt k\neq j}}\frac{1}{z_j-z_k}
          \left[ {\rm Tr}A_{nj}{\rm Tr}A_{nk}-{\rm Tr}A_{nj}A_{nk}-A_{nj,11}-A_{nk,11} \right], \quad j=1,\ldots,N .
\end{equation}
We substitute (\ref{An_resJ}) and (\ref{An_res0}) for the residues of the spectral 
matrix in the above formula in conjunction with the transition relations (\ref{Tform:a}) 
and (\ref{Tform:b}). Into the resulting expression we must further substitute the 
representations of the spectral coefficients (\ref{ThetaRep}), (\ref{OmegaRep}) and 
(\ref{2VRep}). Some care needs to be exercised to ensure that the $ j,k=0 $ contributions 
are correctly accounted for and this leads to a separation between these terms and 
the remaining generic ones. In the group of generic terms there is a summation over $ z_k $ 
for $ 0\leq k\leq N+1 $ and in these terms we can employ the evaluations given in
(\ref{Ssum:a}), (\ref{Ssum:b}), (\ref{Ssum:c}), (\ref{Ssum:d}), (\ref{Ssum:e}) and 
(\ref{Ssum:f}). After reinstating these terms in the total expression considerable 
cancellation is effected. Still further simplification is possible by using the 
transcendent sums (\ref{Tsum:b}), (\ref{Tsum:c}), (\ref{Tsum:d}), (\ref{Tsum:h}) 
and lastly (\ref{Tsum:f}). The final result yields (\ref{Ham}) which is precisely 
Eq.(4.3.7) in \cite{IKSY_1991} after making all the notational correspondences. The 
verification of (\ref{Ham_qDer}) and (\ref{Ham_pDer}) using this Hamiltonian is a 
straight forward calculation.
\end{proof}

\begin{remark}
The Riemann-Papperitz symbol (see subsection 1.4 of \cite{IKSY_1991}, "Fuchsian Equations" 
for the definition) for our system is
\begin{equation}
\left\{
\begin{array}{ccccccccc}
	z_0=0    & z_1     & \cdots & z_N     & z_{N+1}=1 & z_{N+2}=\infty             & q_1 & \cdots & q_N \\
	0        &  0      & \cdots & 0       & 0         & -n                         & 0   & \cdots & 0 \\
	n-\rho_0 & -\rho_1 & \cdots & -\rho_N & -\rho     & 1+\sum^{N+1}_{j=0}\rho_j   & 2   & \cdots & 2 
\end{array}
\right\} .
\end{equation}
\end{remark}

\begin{remark}
One could formulate the Hamiltonian system using the "conjugate" set of variables $ \Theta^*_n, \Omega^*_n $
instead, and in fact all of the dynamical description, however we do not carry out this task as no new
understanding would result from it.
\end{remark}

\begin{remark}
This Hamiltonian system is not polynomial in the canonical variables $ q_r $ and the dynamical 
equations for $ q_r $ do not possess the Painlev\'e property in $ z_j $, however using the well-known 
canonical transformation to the new Hamiltonian system $ {\mathcal H}_N = \{Q_j,P_j;H_j,t_j\} $ 
\cite{KO_1984}, \cite{IKSY_1991}
\begin{align}
   t_j & = \frac{z_j}{z_j-1} ,
\label{Hfix:a} \\
   Q_j & = \frac{z_j\Theta_n(z_j)}{\Theta_{\infty}W'(z_j)} ,
\label{Hfix:b} \\
   P_j & = -(z_j-1)\sum^{N}_{r=1}\frac{p_r}{q_r(q_r-1)}
    \frac{\Theta_{\infty}W(q_r)}{(q_r-z_j)\Theta'_n(q_r)} ,
\label{Hfix:c}
\end{align}
both these deficiencies can be removed.
\end{remark}

\begin{remark}
The definitions of the co-ordinates (\ref{ThetaRep}) and (\ref{Hamp}), (\ref{Ham}) bear 
some resemblance to those in Syklanin's separation of variables procedure \cite{HW_1995a}, \cite{HW_1995b}
although we cannot offer any explanation of this here.
\end{remark}

\section{The discrete Garnier Equations for the $ L(1^{M+1};2) $ Garnier Systems}\label{SemiClassical}
\setcounter{equation}{0}

In this section we derive recurrences in $ n $ or equivalently when 
$ \theta_0 \mapsto \theta_0 \pm 1 $ and $ \theta_{\infty} \mapsto \theta_{\infty} \pm 1 $
for the appropriate variables, thereby characterising the system in this
way. We will find the discrete fifth Painlev\'e equation as the simplest case
corresponding to one free deformation variable and the higher analogues of this
recurrence system for the multi-variable generalisations. 
These are the discrete Garnier equations.

The system with $ M=3, N=1 $ singularities at the standard positions
\begin{equation}
\left\{
\begin{array}{cccc}
	0 & t & 1 & \infty \\
	n-\rho_0 & -\rho_t & -\rho_1 & n+\rho_0+\rho_t+\rho_1 
\end{array}
\right\} ,
\end{equation}
corresponds to sixth Painlev\'e isomonodromic system and was treated extensively 
in \cite{FW_2004b}, especially with regard to the forms of the discrete Painlev\'e 
equations. The weight data is
\begin{equation}
   W(z) = z(z-t)(z-1) = z^3-e_1z^2+e_2z-e_3, \quad
  2V(z) = W\sum_{j=0,t,1}\frac{\rho_j}{z-z_j} = m_0z^2-m_1z+m_2 .
\end{equation}
The spectral coefficients can be parameterised in the following way
\begin{gather}
  \frac{\kappa_{n+1}}{\kappa_n}\Theta_n(z) = \vartheta_n + (n+1+m_0)z ,
  \\
  \Omega_n(z) = - ne_2+\frac{1}{2}m_2 - (\omega_n-\frac{1}{2}m_1)z + (1+\frac{1}{2}m_0)z^2 ,
\end{gather}
whilst the sub-leading coefficients can be related to the polynomial coefficients 
themselves by 
\begin{gather}
  \vartheta_n = -\frac{r_{n}}{r_{n+1}}(n-\rho_0)t ,
  \\
  \omega_n = 1+\rho_0+\rho_t+(1+\rho_0+\rho_1)t+(n+2+m_0)\left[\lambda_{n+2}-\frac{r_{n+2}}{r_{n+1}}\right]
                          -(n+1+m_0)\lambda_{n+1} .
\end{gather}

\begin{proposition}[\cite{FW_2004b}]\label{discreteG3}
The $ n $-recurrence for the bi-orthogonal polynomial system with the $ M=3 $
regular semi-classical weight is governed by the system of coupled first order 
difference equations
\begin{equation}                  
  tf_nf_{n+1}
  = \frac{[\omega_n+n-t-\rho_0(t+1)-(\rho_t+\rho_1)t][\omega_n+n-t-\rho_0(t+1)-\rho_t-\rho_1t]}
           {[\omega_n+nt-1-\rho_0(t+1)-\rho_t-\rho_1][\omega_n+nt-1-\rho_0(t+1)-\rho_t-\rho_1t]} ,
\end{equation}
and
\begin{equation}       
  \omega_n+\omega_{n-1}+(2n-1)t-2-2\rho_0(t+1)-2\rho_t-\rho_1(t+1) 
  = (n-\rho_0)\frac{1-t}{f_n-1}+(n+1+m_0)\frac{1-t}{tf_n-1} .
\end{equation}
The transformations relating these variables to the bi-orthogonal system are
given by
\begin{equation}
  tf_n := \frac{\Theta_n(t)}{\Theta_n(1)} .
\end{equation}
\end{proposition}
\begin{proof}
We refer the reader to the proof of the next case, Proposition \ref{discreteG4}, as
the methods employed in both cases are the same.
\end{proof}

The above coupled recurrence system is equivalent to the canonical "discrete fifth 
Painlev\'e equation" \cite{GORS_1998}, \cite{NRGO_2001} with the mapping
\begin{equation}
   t \mapsto 1/t, \quad
   \omega_n \mapsto (1-t)\omega_n-nt+1+\rho_0(t+1)+\rho_t+\rho_1 ,
\end{equation}
and the identification
\begin{equation}
   \alpha_0 = \rho_t, \quad
   \alpha_1 = n-\rho_0, \quad
   \alpha_2 = -n-\rho_t-\rho_1, \quad
   \alpha_3 = \rho_1, \quad
   \alpha_4 = n+1+\rho_0+\rho_t+\rho_1 .
\end{equation}

In the case $ M=4, N=2 $ we have the two-variable generalisation of the sixth Painlev\'e
equation or the $ L(1^5;2) $ two-variable Garnier system. A standard placement of
the singularities, without loss of generality, would be
\begin{equation}
\left\{
\begin{array}{ccccc}
    0 & s &  t & 1 & \infty \\
   n-\rho_0 & -\rho_s & -\rho_t & -\rho_1 & 
   n+\rho_0+\rho_s+\rho_t+\rho_1
\end{array}
\right\} .
\end{equation} 
We write for notational convenience the weight data in the following way
\begin{gather}
   W(z) = z(z-1)(z-s)(z-t) = z^4-e_1z^3+e_2z^2-e_3z+e_4 ,
\label{m=4_W} \\
  2V(z) = W\sum_{j=0,s,t,1}\frac{\rho_j}{z-z_j} = m_0z^3-m_1z^2+m_2z-m_3 .
\label{m=4_2V}
\end{gather}
The Toeplitz elements satisfy the third order linear difference equation
\begin{multline}
  (j-\rho_0)st w_{j} 
   -\left[ (j-1-\rho_0)(s+t+st)-\rho_1 st-\rho_s t-\rho_t s \right]w_{j-1}
 \\
   +\left[ (j-2-\rho_0)(1+s+t)-\rho_1(s+t)-\rho_s(t+1)-\rho_t(s+1) \right]w_{j-2}
 \\
   -(j-3-\rho_0-\rho_1-\rho_s-\rho_t)w_{j-3} = 0 ,\quad j\in \Z ,
\end{multline}
and we take $ w_{-1},w_0,w_1 $ as defining a solution of the system. 
Also, according to (\ref{Thexp:a})-(\ref{Omexp:b}), we can parameterise the spectral
coefficients as 
\begin{gather}
  \frac{\kappa_{n+1}}{\kappa_n}\Theta_n(z) =
     (n-\rho_0)e_3\frac{r_n}{r_{n+1}}+\vartheta_n z+(n+1+m_0) z^2 ,
\label{L2_Theta:a} \\
  \Omega_n(z) = ne_3-\frac{1}{2}m_3 + (\omega_n-\frac{1}{2}m_2) z
   + (\varpi_n+\frac{1}{2}m_1) z^2 + (1+\frac{1}{2}m_0) z^3 ,
\label{L2_Omega:a}
\end{gather}
in terms of $ \vartheta_n $, $ \omega_n $ and $ \varpi_n $.

We can relate the new variables introduced above to the coefficients of the polynomials
through the explicit forms of the spectral coefficients (\ref{Thexp:a})-(\ref{Omexp:b})
in the following way
\begin{gather}
   \vartheta_n =
   -(n+1)e_1-m_1+(n+2+m_0)\left[ \frac{r_{n+2}}{r_{n+1}}-r_{n+2}\bar{r}_{n+1} \right]
   -(n+m_0)r_{n+1}\bar{r}_{n}-2\lambda_{n+1} ,
\label{L2_Theta:b} \\
   \omega_n =
    -ne_2+m_2+(n-\rho_0)e_3\left[ \frac{r_n}{r_{n+1}}-\bar{r}_{n+1}r_n \right]-e_3\bar{\lambda}_{n+1} ,
\label{L2_Omega:b} \\
   \varpi_n =
    -e_1-m_1+(n+2+m_0)\left[ \frac{r_{n+2}}{r_{n+1}}-r_{n+2}\bar{r}_{n+1} \right]-\lambda_{n+1} .
\label{L2_Omega:c}
\end{gather}
These formulae only serve to allow the recovery of the original variables and do 
not feature in the recurrence relations.

At this point it is possible to use the foregoing results to derive a system of 
recurrence relations which is the analogue of the fifth discrete Painlev\'e
equation for the two-variable Garnier system.
\begin{proposition}\label{discreteG4}
The following system of coupled first order recurrence relations 
in $ n $ for the variables 
$ \{ f_{n},g_{n},\omega_n,\varpi_n \}^{\infty}_{n=0} $ completely 
characterises the bi-orthogonal polynomial system
\begin{equation} 
  sf_{n}f_{n+1} =
  \frac{\left[ \omega_n+s\varpi_n+(1+m_0)s^2+(n-\rho_0)t+\rho_s(s-t)(1-s) \right]
        \left[ \omega_n+s\varpi_n+(1+m_0)s^2+(n-\rho_0)t \right]}
       {\left[ \omega_n+\varpi_n+1+m_0+(n-\rho_0)st+\rho_1(1-s)(t-1) \right]
        \left[ \omega_n+\varpi_n+1+m_0+(n-\rho_0)st \right]} ,
\label{L2_dP:a}
\end{equation} 
\begin{equation} 
  tg_{n}g_{n+1} =
  \frac{\left[ \omega_n+t\varpi_n+(1+m_0)t^2+(n-\rho_0)s+\rho_t(t-s)(1-t) \right]
        \left[ \omega_n+t\varpi_n+(1+m_0)t^2+(n-\rho_0)s \right]}
       {\left[ \omega_n+\varpi_n+1+m_0+(n-\rho_0)st+\rho_1(1-s)(t-1) \right]
        \left[ \omega_n+\varpi_n+1+m_0+(n-\rho_0)st \right]} ,
\label{L2_dP:b}
\end{equation} 
\begin{multline} 
  \omega_{n}+\omega_{n-1} = m_2-(n-1)(s+t+st) \\
  + (n-\rho_0)
    \frac{s^2-t^2+(1-s^2)tg_{n}-(1-t^2)sf_{n}}{t-s+(1-t)f_{n}-(1-s)g_{n}}
  + (n+1+m_0)st
    \frac{t-s+(1-t)f_{n}-(1-s)g_{n}}{t-s+(1-t)sf_{n}-(1-s)tg_{n}} ,
\label{L2_dP:c}
\end{multline} 
\begin{multline} 
  \varpi_{n}+\varpi_{n-1} = -m_1+(n-1)(1+s+t) \\
  + (n-\rho_0)
    \frac{t-s+(1-t)sf_{n}-(1-s)tg_{n}}{t-s+(1-t)f_{n}-(1-s)g_{n}}
  + (n+1+m_0)
    \frac{s^2-t^2+(1-s^2)tg_{n}-(1-t^2)sf_{n}}{t-s+(1-t)sf_{n}-(1-s)tg_{n}} ,
\label{L2_dP:d}
\end{multline} 
where
\begin{equation}
   sf_{n} := \frac{\Theta_{n}(s)}{\Theta_{n}(1)}, \quad
   tg_{n} := \frac{\Theta_{n}(t)}{\Theta_{n}(1)}.
\label{L2_fgDefn}
\end{equation}
The recurrence relations are subject to the initial conditions
\begin{align}
  f_0 & = 
  \frac{(1+m_0)sw_{-1}-\left[ \rho_0 s(t+1)+\rho_1 st+\rho_s t+\rho_t s \right]w_0-(1-\rho_0)stw_{1}}
       {(1+m_0)w_{-1}-\left[ \rho_0 (s+t)+\rho_1 st+\rho_s t+\rho_t s \right]w_0-(1-\rho_0)stw_{1}} ,
  \\
  g_0 & =
  \frac{(1+m_0)tw_{-1}-\left[ \rho_0 t(s+1)+\rho_1 st+\rho_s t+\rho_t s \right]w_0-(1-\rho_0)stw_{1}}
       {(1+m_0)w_{-1}-\left[ \rho_0 (s+t)+\rho_1 st+\rho_s t+\rho_t s \right]w_0-(1-\rho_0)stw_{1}} ,
  \\
  \omega_0 & = 
  (1-\rho_0)st\frac{w_{1}}{w_{0}}+\rho_0 st\frac{w_{0}}{w_{-1}}
        +\rho_0(s+t+st)+\rho_1 st+\rho_s t+\rho_t s ,
  \\
  \varpi_0 & = 
  -(1+m_0)\frac{w_{-1}}{w_{0}}-(1-\rho_0)st\frac{w_{1}}{w_{-1}}
        -\left[ \rho_0(s+t+st)+\rho_1 st+\rho_s t+\rho_t s \right]\frac{w_{0}}{w_{-1}} .
\end{align}
\end{proposition}
\begin{proof}
The first two relations (\ref{L2_dP:a},\ref{L2_dP:b}) follow from the 
evaluation of (\ref{OTeq:a}) at $ j=s,t,1 $ and taking ratios. From
the definitions (\ref{L2_fgDefn}) and (\ref{L2_Theta:a}) we note that
\begin{align} 
  f_{n} & =
  \frac{(n-\rho_0)t\frac{\DySt r_n}{\DySt r_{n+1}}+\vartheta_n+(n+1+m_0)s} 
       {(n-\rho_0)st\frac{\DySt r_n}{\DySt r_{n+1}}+\vartheta_n+n+1+m_0} ,
  \\ 
  g_{n} & =
  \frac{(n-\rho_0)s\frac{\DySt r_n}{\DySt r_{n+1}}+\vartheta_n+(n+1+m_0)t} 
       {(n-\rho_0)st\frac{\DySt r_n}{\DySt r_{n+1}}+\vartheta_n+n+1+m_0} .
\end{align} 
Solving for $ \vartheta_n $ and $ r_n/r_{n+1} $ in terms of $ f_n,g_n $
we find that
\begin{equation}
   \vartheta_n = (n+1+m_0)
    \frac{s^2-t^2+(1-s^2)tg_{n}-(1-t^2)sf_{n}}{t-s+(1-t)sf_{n}-(1-s)tg_{n}} ,
\end{equation}
and
\begin{equation}
   \frac{r_n}{r_{n+1}} = \frac{n+1+m_0}{n-\rho_0}
    \frac{t-s+(1-t)f_{n}-(1-s)g_{n}}{t-s+(1-t)sf_{n}-(1-s)tg_{n}} .
\end{equation}
The second pair of relations (\ref{L2_dP:c},\ref{L2_dP:d}) follow from
the evaluation of the recurrence (\ref{rrCf:a}) at $ z=s,t $, then
solving for their left-hand sides and utilising the preceding expressions.
\end{proof}

We now give the result for an arbitrary number of variables in a canonical placement
expressed by the abbreviated symbol
\begin{equation}
\left\{
\begin{array}{cccccc}
	0 & t_1 & \cdots & t_N & 1 & \infty \\
	n-\rho_0 & -\rho_1 & \cdots & -\rho_N & -\rho & n+\sum^{N+1}_{j=0}\rho_j
\end{array}
\right\} ,
\end{equation}
with $ M=N+2 $. 
To begin with we denote the $ l $-th elementary symmetric function
in the variables $ t_1, \ldots, t_n $ by $ e_l(t_1,\ldots,t_n) $ and the convention
$ e_0=1 $. We adopt the definition and notations for $ W $ and $ 2V $ as given by
(\ref{Wdefn}) and (\ref{2Vdefn}) respectively, along with this renaming of the 
independent variables.
Let us define the set $ T:=\{t_1,\ldots,t_N\} $ and the set omitting the
variable $ t_j $ by $ T_j:=T\backslash \{t_j\} $. Furthermore we define the 
Vandermonde determinant
\begin{equation} 
  \Delta(T) := \prod_{1\leq j<k\leq N}(t_k-t_j) .
\end{equation} 

The defining relation for the weight (\ref{scwgt2}) implies that the
Toeplitz elements $ \{w_k\}^{\infty}_{k=-\infty} $ satisfy the $ N+1 $-order linear 
difference equation
\begin{equation}
  \sum^{N+1}_{l=0} (-)^l\left[ (j-l)e_{N+1-l}-m_{N+1-l} \right]w_{j-l} = 0,
  \quad j\in \Z ,
\label{}
\end{equation} 
with $ N+1 $ consecutive elements being arbitrary "initial" values. These 
"initial" values also define the $ U $ polynomial, which will in turn fix the
initial values of our recurrence relations \ref{discreteGM}, by the expression
\begin{equation}
   U:=\sum^{N+1}_{l=0}u_lz^l ,
\label{}
\end{equation}
where the coefficients are 
\begin{gather}
  u_0 = (-)^Nw_0m_{N+1}, \quad
  u_{N+1} = w_0m_0,
  \\
  u_j = (-)^{N-j}w_0m_{N+1-j}+2\sum^{j-1}_{l=0}(-)^{N+1-l}\left[ (j-l)e_{N+1-l}-m_{N+1-l} \right]w_{j-l},
  \quad j=1,\ldots,N .
\label{}
\end{gather} 

We adopt a parameterisation of the spectral coefficients of the form
\begin{equation}
 \frac{\kappa_{n+1}}{\kappa_{n}}\Theta_n(z)
 = (-)^N(ne_{N+1}-m_{N+1})\frac{r_n}{r_{n+1}}
   +\sum^{N-1}_{l=1}\vartheta^l_n z^l + (n+1+m_0)z^N , 
\label{Thparam}
\end{equation}
and
\begin{equation}
  \Omega_n(z) = (-)^N(ne_{N+1}-\frac{1}{2}m_{N+1})
  +(-)^N\sum^{N}_{l=1}(\omega^l_n+\frac{1}{2}(-)^lm_{N+1-l})z^l + (1+\frac{1}{2}m_0)z^{N+1} ,
\label{Omparam}
\end{equation}
which introduces the variables $ \vartheta^1_n, \ldots, \vartheta^{N-1}_n $ and 
$ \omega^1_n, \ldots, \omega^N_n $.

\begin{proposition}\label{discreteGM}
Define the variables
\begin{equation}
   t_jf^j_{n} := \frac{\Theta_{n}(t_j)}{\Theta_{n}(1)}, \quad j=1,\ldots, N .
\label{f_Defn}
\end{equation}
The following system of $ 2N $ coupled first order recurrence relations 
in $ n $ for the variables $ \{ f^j_n,\omega^j_n \}^{N}_{j=1} $ 
\begin{multline}
  t_jf^j_{n}f^j_{n+1}
	= \frac{\left[
	        (n-\rho_0)\prod_{k\neq j}t_k+\sum^N_{l=1}t^{l-1}_j\omega^l_n+(-)^N(1+m_0)t^N_j
	        \right]}
	        {\left[
	        (n-\rho_0)\prod_{k}t_k+\sum^N_{l=1}\omega^l_n+(-)^N(1+m_0)
	        \right]}
	\\ \times
	  \frac{\left[
	        n\prod_{k\neq j}t_k+\sum^N_{l=1}t^{l-1}_j(\omega^l_n+(-)^lm_{N+1-l})+(-)^Nt^N_j
	        \right]}
	       {\left[
	        n\prod_{k}t_k+\sum^N_{l=1}(\omega^l_n+(-)^lm_{N+1-l})+(-)^N
	        \right]} ,
  \quad j=1,\ldots, N ,
\label{dGarnier:a}
\end{multline} 
and
\begin{multline} 
  \omega^j_{n}+\omega^j_{n-1}+(-)^jm_{N+1-j} = (-)^j(n-1)e_{N+1-j}(T\cup\{1\})
  \\
  + (-)^j(n-\rho_0)
    \frac{\Delta(T)e_{N-j}(T)+\sum^{N}_{l=1}(-)^{N-1+l}t_l\Delta(T_l\cup\{1\})e_{N-j}(T_l\cup\{1\})f^l_n}
         {\Delta(T)+\sum^{N}_{l=1}(-)^{N-1+l}\Delta(T_l\cup\{1\})f^l_n}
  \\
  - (-)^j(n+1+m_0)
    \frac{\Delta(T)e_{N+1-j}(T)+\sum^{N}_{l=1}(-)^{N-1+l}t_l\Delta(T_l\cup\{1\})e_{N+1-j}(T_l\cup\{1\})f^l_n}
         {\Delta(T)+\sum^{N}_{l=1}(-)^{N-1+l}t_l\Delta(T_l\cup\{1\})f^l_n} ,
  \quad j=1,\ldots, N ,
\label{dGarnier:b}
\end{multline} 
completely characterises the bi-orthogonal polynomial system.
The recurrence relations are subject to the initial conditions
\begin{align}
  f^j_0 & = \frac{2V(t_j)-\kappa^2_0U(t_j)}{t_j[2V(1)-\kappa^2_0U(1)]} ,
\label{dGarnier:c}  \\
  (-)^N\omega^j_0 & = -\frac{1}{2}[z^j]\left( 2V+\kappa^2_0U \right) - \frac{w_0}{2w_{-1}}[z^{j-1}]\left( 2V-\kappa^2_0U \right) ,
\label{dGarnier:d} 
\end{align}
for $ j=1,\ldots,N $, where $ [z^l](.) $ denotes the coefficient of $ z^l $ in 
the polynomial.
\end{proposition}
\begin{proof}
To establish (\ref{dGarnier:a}) we utilise (\ref{OTeq:a}) in the following form
\begin{equation}
   \frac{t_j\Theta_n(t_j)\Theta_{n+1}(t_j)}{\Theta_n(1)\Theta_{n+1}(1)}
  = \frac{[\Omega_n(t_j)+V(t_j)][\Omega_n(t_j)-V(t_j)]}{[\Omega_n(1)+V(1)][\Omega_n(1)-V(1)]},
    \quad j=1,\ldots,N.
\end{equation}
We note that the factors appearing on the right-hand side can be found from
\begin{gather}
  \Omega_n+V 
  = (-)^N(n-\rho_0)e_{N+1}+(-)^N\sum^N_{l=1}\omega^l_n z^l+(1+m_0)z^{N+1} ,
  \\
  \Omega_n-V 
  = (-)^Nne_{N+1}+(-)^N\sum^N_{l=1}(\omega^l_n+(-)^lm_{N+1-l}) z^l+z^{N+1} .
\end{gather}
Then (\ref{dGarnier:a}) follows from this result and the definition (\ref{f_Defn}).

To prove the relations (\ref{dGarnier:b}) we first need to invert the definition 
(\ref{f_Defn}) along with the parameterisation (\ref{Thparam}) for the coefficients 
$ r_{n}/r_{n+1} $ and $ \vartheta^j_n $ for $ j=1,\ldots,N $. We find
\begin{equation} 
  (n-\rho_0)\frac{r_{n}}{r_{n+1}}
  = (n+1+m_0)
    \frac{\Delta(T)+\sum^{N}_{l=1}(-)^{N-1+l}\Delta(T_l\cup\{1\})f^l_n}
         {\Delta(T)+\sum^{N}_{l=1}(-)^{N-1+l}t_l\Delta(T_l\cup\{1\})f^l_n} ,
\label{dGaux:a}
\end{equation}
and
\begin{equation} 
  \vartheta^j_{n}
  = (-)^{N+j}(n+1+m_0)
    \frac{\Delta(T)e_{N-j}(T)+\sum^{N}_{l=1}(-)^{N-1+l}t_l\Delta(T_l\cup\{1\})e_{N-j}(T_l\cup\{1\})f^l_n}
         {\Delta(T)+\sum^{N}_{l=1}(-)^{N-1+l}t_l\Delta(T_l\cup\{1\})f^l_n} ,
\label{dGaux:b}
\end{equation}
for $ j=1,\ldots,N-1 $
by using the identity for Vandermonde determinants $ \Delta_j(T) $ with the $ j $-th 
column missing
\begin{equation}
  \Delta_j(T) := \det
       \begin{pmatrix}
              1 & t_1 & \cdots & [] & \cdots & t^N_1 \cr
              \vdots & \vdots & & \vdots & & \vdots \cr
              1 & t_N & \cdots & [] & \cdots & t^N_N \cr
       \end{pmatrix}
               = e_{N-j}(T)\Delta_N(T) ,
\end{equation} 
for $ j=0,\ldots,N $.
The initial values of the dependent variables (\ref{dGarnier:c} and (\ref{dGarnier:d}) 
are found using (\ref{Theta0}) and (\ref{Omega0}) respectively and their definitions.
\end{proof}

\begin{remark}
Of primary interest in applications are the Toeplitz determinants (\ref{Uavge}) which are
also $ \tau$-functions in the sense of Jimbo-Miwa-Ueno's definition 
\cite{JM_1981a},\cite{JM_1981b},\cite{JM_1981c}.
Recovery of the $ \tau$-function $ I_n $ from the $ f^j_n, \omega^j_n $ of the previous proposition 
can be achieved by solving the recurrence in $ r_n $ using (\ref{dGaux:a}) and then by
employing either of
\begin{equation}
  \vartheta^{N-1}_{n} = -(n+1)e_1-m_1
   +(n+2+m_0)\left( \frac{r_{n+2}}{r_{n+1}}-\lambda_{n+2} \right)+(n+m_0)\lambda_n ,
\end{equation}
or
\begin{equation}
  (-)^N\omega^{N}_n = -e_1-m_1
   +(n+1+m_0)\lambda_{n+1}+(n+2+m_0)\left( \frac{r_{n+2}}{r_{n+1}}-\lambda_{n+2} \right) ,
\end{equation}
and solving these recurrences for $ \lambda_n $. This latter sequence can then used to
find $ \bar{r}_n $ via (\ref{l}), and with both of the $r$-coefficients one can
use (\ref{I0}) to solve for the recurrence in $ I_n $. 
\end{remark}

\begin{remark}
The mapping between the Hamiltonian variables for $ {\mathcal G}_N $ and those of the above proposition 
are given by (\ref{f_Defn}) with (\ref{ThetaRep}) and
\begin{equation}
   -W(q_r)p_r = (-)^N(n-\rho_0)e_{N+1}+(-1)^N\sum^{N}_{j=1}\omega^j_n q^j_r+(1+m_0)q^{N+1}_r ,
\end{equation}
or
\begin{equation}
   (-)^j\omega^j_n 
  = \left[ 1+m_0+(n-\rho_0)\frac{e_N(T)}{e_N({\mathcal Q})} \right]e_{N-j}({\mathcal Q})
    -\sum^{N}_{r=1} e_{N-j}({\mathcal Q}_r)\frac{(q_r-1)\prod^N_{k=1}(q_r-t_k)}{\prod_{s\neq r}(q_r-q_s)}p_r ,
\end{equation}
where $ {\mathcal Q}:=\{q_s\}^N_{s=1} $ and $ {\mathcal Q}_r:={\mathcal Q}\backslash \{q_r\} $.
Written in these Hamiltonian co-ordinates the second of the coupled recurrences 
(\ref{dGarnier:b}) is
\begin{equation}
   p_{n+1,r}+p_{n,r} = \frac{n}{q_r}-\frac{2V(q_r)}{W(q_r)} .
\end{equation}
\end{remark}

\begin{remark}
A system of discrete Garnier equations could be formulated using the "conjugate"
variables $ \Theta^*_n $ and $ \Omega^*_n $ but this would not differ in essence from 
the one presented here.
\end{remark}

\begin{remark}
We have not attempted to check whether singularity confinement \cite{GNR_1999}, 
an algebraic entropy criteria \cite{HV_1998}, \cite{BV_1999}, \cite{TaT_2001} or one based on 
Nevalinna theory \cite{AHH_2000}, \cite{Ha_2005} applies to our recurrences and this remains an 
outstanding issue.
\end{remark}

\begin{remark}
There have been reports of recurrence relations for the Garnier systems in 
\cite{TN_2006} and \cite{NW_2001}, however the explicit relationship between these 
equations and the recurrences reported here remains to be elucidated.
\end{remark}

\begin{remark}
Another possible method for deriving the recurrence relations of Proposition \ref{discreteGM}
would be to construct the Schlesinger transformation, or lattice translation, operators
from the fundamental reflection and automorphism operators of the affine Weyl group
$ B^{(1)}_{N+3} $, which have been studied in \cite{Ki_1990}, \cite{Ts_2003a}, \cite{Su_2005}.
This task was carried out for $ M=3, N=1 $, i.e. for the sixth Painlev\'e equation in
\cite{FW_2006a}, however this involved a laborious calculation and it may not be feasible to
employ this approach to the multi-variable extension.
\end{remark}

In this study we have focused on a particular type of transformation, namely that
of $ n \mapsto n\pm 1 $ or equivalently $ \theta_0 \mapsto \theta_0\pm 1 $ and
$ \theta_{\infty} \mapsto \theta_{\infty}\pm 1 $, which is natural within this 
context. However one could legitimately ask for the recurrence systems for the
transformations $ \rho_j \mapsto \rho_j\pm 1 $, which are part of the larger group
of Schlesinger transformations. The theory of these, in the context of bi-orthogonal
system on the unit circle, has been investigated in \cite{Wi_2009b} but analogues of the
recurrences found here were not given there.

\section{Acknowledgments}
This research has been supported by the Australian Research Council (ARC) and partially 
supported by the ARC Centre of Excellence for Mathematics and Statistics of Complex 
Systems. The author also appreciates the assistance and advice given by Chris Cosgrove
regarding the implementation of the recurrence relations in computer algebra systems. 

\setcounter{equation}{0}
\bibliographystyle{amsplain}
\bibliography{moment,nonlinear,random_matrices}

\end{document}